\documentclass{article}
\usepackage{amsmath}
\usepackage{graphicx}
\usepackage{amsfonts}
\usepackage{amssymb}
\newtheorem{theorem}{Theorem}[section]

\newtheorem{corollary}[theorem]{Corollary}

\newtheorem{lemma}[theorem]{Lemma}

\newtheorem{proposition}[theorem]{Proposition}

\newenvironment{proof}[1][Proof]{\textbf{#1.} }{\ \rule{0.5em}{0.5em}}
\newdimen\dummy
\dummy=\oddsidemargin
\begin{document}

\title{Explicit Hilbert Spaces for certain unipotent representations III}
\author{Alexander Dvorsky\thanks{\texttt{dvorsky@math.miami.edu, }Department of
Mathematics, University of Miami, Coral Gables, FL 33126}\ and Siddhartha
Sahi\thanks{\texttt{sahi@math.rutgers.edu}, Department of Mathematics, Rutgers
University, New Brunswick, NJ 08903}}
\maketitle

\section{Introduction}

This paper is the culmination of a series dedicated to the problem of
constructing explicit analytic models for small unitary representations of
certain semisimple Lie groups. (See \cite{sahi-expl}, \cite{shilov},
\cite{sahi-dp}, \cite{tens}, \cite{hilbert-two} and also \cite{sahi-stein},
\cite{kostant-sahi}.)

The groups $G$ that we consider are those that arise from real semisimple
Jordan algebras via the Tits-Koecher-Kantor construction. Such a $G$ is
characterized by the fact that it admits a parabolic subgroup $P=LN$ which is
conjugate to its opposite $\overline{P}=L\overline{N}$, and for which the
nilradicals $N$ and $\overline{N}$ are abelian. We identify $N$ and
$\overline{N}$ with their Lie algebras $\frak{n}$ and $\overline{\frak{n}}$
via the exponential map, and moreover we identify $\overline{\frak{n}}$ with
$\frak{n}^{\ast}$ via a suitably rescaled Killing form. The space
$\overline{\frak{n}}$ is naturally endowed with a real semisimple Jordan
algebra structure, which is isotopic to the one giving rise to $G$ and $P$. We
write $n$ for the rank of $\overline{\frak{n}}$ as a Jordan algebra.

In this situation, the Levi component $L$ has a finite number of orbits on
$\overline{\frak{n}}$; each orbit has a \emph{rank } $\leq n$, and carries a
measure $d\mu$ which transforms by a character under $L$. For the non-open
orbits (rank $<n$) this measure is unique up to scalar multiple; for the open
orbits (rank $=n$) there is a one-parameter family of such measures. For each
non-open orbit $\mathcal{O}$ we consider the Hilbert space $\mathcal{H}%
_{\mathcal{O}}=\mathcal{L}^{2}(\mathcal{O},d\mu)$. By\ Mackey theory this
space carries a natural irreducible unitary representation $\pi_{\mathcal{O}}$
of $P$, and we consider the following two problems:

\begin{itemize}
\item  Extend $\pi_{\mathcal{O}}$ to a unitary representation of $G$.

\item  Decompose $\pi_{\mathcal{O}^{1}}\otimes\cdots\otimes\pi_{\mathcal{O}s}
$, for $\operatorname*{rank}\mathcal{O}^{1}+...+\operatorname*{rank}%
\mathcal{O}^{s}\leq n $.
\end{itemize}

For Euclidean Jordan algebras the first problem was solved in \cite{sahi-expl}
and \cite{shilov}, and the second problem was solved in \cite{tens}. Therefore
in this paper we restrict our attention to non-Euclidean Jordan algebras. The
case of rank $1$ orbits in non-Euclidean Jordan algebras was completely
settled in \cite{hilbert-two}. It turned out that for rank $1$ orbits in
certain rank $2$ Jordan algebras the representation $\pi_{\mathcal{O}}$ cannot
be extended to $G$. We shall call these orbits \emph{inadmissible}; they arise
for the groups $G=O(p,q),\;p\neq q$ and are connected with the Howe-Vogan
``no-go'' result for minimal representations of such groups (\cite{vogan-nogo}%
). For rank $1$ orbits in the remaining Jordan algebras we established the
desired extension and calculated the tensor product decomposition in
\cite{hilbert-two}.

In this paper we consider the general case, and prove the following results:

First suppose $\mathcal{O}$ is an admissible non-open orbit in a non-Euclidean
Jordan algebra $\overline{N}$.

\begin{theorem}
For $\mathcal{O}$ as above, the representation $\pi_{\mathcal{O}}$ of $P$
extends to an irreducible spherical unitary representation of $G$ on
$\mathcal{H}_{\mathcal{O}}$.\label{theoremA}
\end{theorem}

Now suppose $\mathcal{O}^{1},\ldots,\mathcal{O}^{s}$ are admissible non-open
orbits in $\overline{N}$ such that $\operatorname*{rank}\mathcal{O}%
^{1}+...+\operatorname*{rank}\mathcal{O}^{s}\leq n$. In subsection
\ref{sec-tensor} we define a reductive homogeneous space $G^{\prime}%
/H^{\prime}$, essentially the generic fiber of the addition map from
$\mathcal{O}^{1}\times\ldots\times\mathcal{O}^{s}$ to $\overline{N}$, and
consider the decomposition of the quasi-regular representation
\[
\mathcal{L}^{2}(G^{\prime}/H^{\prime})=\int_{\widehat{G^{\prime}}}^{\oplus
}m(\sigma)\sigma\,d\rho(\sigma),
\]
where $m(\sigma)$ is the multiplicity function and $d\rho(\sigma)$ is the
Plancherel measure.

\begin{theorem}
Let $\mathcal{O}^{1},\ldots,\mathcal{O}^{s}$ and $G^{\prime},H^{\prime}$ be as
above; then there is a map $\mathbf{\theta}$ from the $H^{\prime}$-spherical
dual of $G^{\prime}$ to the unitary dual of $G$ such that
\[
\pi_{\mathcal{O}^{1}}\otimes\cdots\otimes\pi_{\mathcal{O}^{s}}=\int
_{\widehat{G^{\prime}}}^{\oplus}m(\sigma)\mathbf{\theta}(\sigma)\,d\rho
(\sigma)\text{.\label{theoremB}}%
\]
\end{theorem}

Our approach to these results requires consideration of three different kinds
of groups, each with its own flavor of representation theory. These are:

\begin{itemize}
\item  Harish-Chandra modules for semisimple Lie groups.

\item  Operator algebras for parabolic subgroups.

\item  Fourier analysis for abelian nilradicals.
\end{itemize}

The Harish-Chandra theory was studied in \cite{sahi-dp}. The necessary
operator-algebraic results ($C^{\ast}$-algebras, von Neumann algebras) were
already obtained in \cite{tens} and \cite{hilbert-two}. Thus the missing
ingredient, which is provided by this paper, involves abelian Fourier
analysis. The key result, Proposition \ref{=prop-L2}, shows that a certain
function (eventually, the ``spherical'' vector in $\pi_{\mathcal{O}}$) belongs
to $\mathcal{L}^{2}(\mathcal{O},d\mu)$. For rank $1$ orbits this result was
obtained in \cite{hilbert-two} by establishing a close connection between this
function and a certain one-variable Bessel $K$-function. The required
$\mathcal{L}^{2}$-estimate was then deduced from a precise knowledge of the
singularity of the Bessel $K$-function at $0$.

For higher rank orbits, we expect that there should exist a similar connection
between the spherical vectors and multivariate Bessel $K$-functions. However
in order to exploit this connection one would have to first develop the theory
of such functions, possibly along the lines of the theory of the multivariate
Bessel $J$-functions of \cite{opdam}. While we feel that the connection with
multivariate Bessel $K$-functions is of interest and should be pursued
further, in the present paper we follow a different approach that allows us to
obtain the desired estimate directly, obviating the need to first study Bessel
functions. The key here is a ``stability'' result (Lemma \ref{lemma-stab})
which transfers the problem from a non-open orbit to a related problem on the
open orbit for a smaller group. The open orbit problem turns out to be easier
to solve.

This approach was inspired in part by a recent paper of Shimura \cite{shimura}%
. We wish to thank L. Barchini for drawing our attention to this paper.

\section{Preliminaries}

In this section we recall basic facts about the Tits-Kantor-Koecher
construction. All results of this section are well-known. More details may be
found in \cite{kostant-sahi}, \cite{hilbert-two} and in the references therein
(in particular, \cite{braun} and \cite{loos}).

\subsection{The pair $\left(  G,P\right)  $}

The Tits-Kantor-Koecher construction associates to a real simple Jordan
algebra, a pair $\left(  G,P\right)  $, where $G$ is a real simple Lie group
with Cartan involution $\theta$, and maximal compact subgroup $K$; and $P$ is
a parabolic subgroup with Levi decomposition $LN$, say. In the context of Lie
theory, these pairs can be characterized as follows:

\begin{enumerate}
\item $N$ is abelian.

\item $P$ is $G$-conjugate to its opposite parabolic $\overline{P}%
=\theta(P)=L\overline{N}$.
\end{enumerate}

Conversely, in the above situation one can endow $N$ with the structure of a
real simple Jordan algebra, unique up to the choice of an identity element.

In view of the classification of real simple Lie groups and their parabolic
subgroups in terms of the possible restricted root systems, it is an easy
matter to determine the above pairs. The two conditions correspond to simple
restricted roots $\alpha$ such that

\begin{enumerate}
\item $\alpha$ has coefficient $1$ in the highest root

\item $\alpha$ satisfies $\alpha=-w_{0}\alpha$ for the longest element $w_{0}$
of the Weyl group.
\end{enumerate}

The conditions 1$.$ and 2$.$ each give rise to a \emph{symmetric} space
denoted by $K/M$ and $L/H$ respectively, and much of the relevant information
about the Jordan algebra and the associated pair $\left(  G,P\right)  $ can be
described in a simple and coherent manner in terms of these symmetric spaces.
This makes it possible to have a uniform discussion for the most part, with
only some occasional arguments requiring case-by-case considerations.

In the next few sections we describe these spaces and conclude by giving a
complete list of examples.

We follow the practice of denoting the real Lie algebras of various Lie groups
by the corresponding fraktur letters; with the exception of $\frak{p}$ which
will denote instead the $-1$ eigenspace of $\theta$ in the Cartan
decomposition $\frak{g}=\frak{k}\oplus\frak{p}$.

\subsection{The symmetric space $K/M$}

Condition 1$.$ implies (and is equivalent to the assertion) that the subgroup
$L$ is a \emph{symmetric} subgroup of $G$, and $M=K\cap L$ is a symmetric
subgroup of $K$. \ Let $\frak{t}$ be a maximal toral subalgebra for the
compact symmetric space $K/M$, \textit{i.e.} a Cartan subspace for in the
orthogonal complement of $\frak{m}$ in $\frak{k}$. \ The real rank of $N$ as a
Jordan algebra is $n=\dim_{\mathbb{R}}\frak{t}$.

The roots of $\frak{t}_{\mathbb{C}}$ in $\frak{g}_{\mathbb{C}}$ \emph{always}
form a root system of type $C_{n}$, and we fix a basis $\{\gamma_{1}%
,\gamma_{2},\ldots,\gamma_{n}\}$ of $\frak{t}^{\ast}$ such that
\[
\Sigma(\frak{t}_{\mathbb{C}},\frak{g}_{\mathbb{C}})=\{\pm(\gamma_{i}\pm
\gamma_{j})/2,\pm\gamma_{j}\}\text{.}%
\]
For the subsystem $\Sigma=\Sigma(\frak{t}_{\mathbb{C}},\frak{k}_{\mathbb{C}}%
)$, there are \emph{three }possibilities:
\[
A_{n-1}=\{\pm(\gamma_{i}-\gamma_{j})/2\}\text{, \ }D_{n}=\{\pm(\gamma_{i}%
\pm\gamma_{j})/2\}\text{, and }C_{n}.
\]

The first of these cases arises precisely when $N$ is a \emph{Euclidean}
Jordan algebra. This case was studied in \cite{sahi-expl}, therefore we
restrict our attention to the last two cases. If $\Sigma$ is $C_{n}$, there
are two multiplicities, corresponding to the short and long roots, which we
denote by $d$ and $e,$ respectively. If $\Sigma$ is $D_{n}$, and $n\neq2,$
then there is a single multiplicity, which we denote by $d$, so that $D_{n}$
may be regarded as a special case of $C_{n}$, with $e=0$.

The root system $D_{2}$ $\approx A_{1}\times A_{1}$ is reducible and there are
two (possibly different) root multiplicities. In what follows, we
explicitly\emph{\ }exclude the case when these multiplicities are
different.\emph{\ }This means that we exclude from consideration the groups
\[
G=O(p,q),N=\mathbb{R}^{p-1,q-1}(p\neq q)\text{.}%
\]
\emph{\ }When the two multiplicities coincide $(p=q)$, we once again denote
the common multiplicity by $d$.\label{sec-KM}

\subsection{$S$-triples and the Cayley transform}

The discussion of the various cases can be made uniform by emphasizing the
special role played by a family of $n$ commuting $SL_{2}$'s or $S$-triples,
together with the associated Cayley transform.

For the Lie algebra $\frak{sl}_{2}(\mathbb{C})$, we define
\begin{align*}
x  &  =\left[
\begin{array}
[c]{cc}%
0 & 1\\
0 & 0
\end{array}
\right]  ,\text{ }y=\left[
\begin{array}
[c]{cc}%
0 & 0\\
1 & 0
\end{array}
\right]  ,\text{ }h=\left[
\begin{array}
[c]{cc}%
1 & 0\\
0 & -1
\end{array}
\right]  ,\\
X  &  =\frac{1}{2}\left[
\begin{array}
[c]{cc}%
i & 1\\
1 & -i
\end{array}
\right]  ,Y=\frac{1}{2}\left[
\begin{array}
[c]{cc}%
-i & 1\\
1 & i
\end{array}
\right]  ,\text{ }H=i\left[
\begin{array}
[c]{cc}%
0 & 1\\
-1 & 0
\end{array}
\right]  .
\end{align*}
The Cayley transform is the automorphism of $\frak{sl}_{2}(\mathbb{C})$ given by%

\[
c=\exp\operatorname*{ad}\frac{\pi i}{4}\left(  X+Y\right)  =\exp
\operatorname*{ad}\frac{\pi i}{4}\left(  x+y\right)  =\text{ }%
\operatorname*{Ad}\frac{1}{\sqrt{2}}\left[
\begin{array}
[c]{cc}%
1 & i\\
i & 1
\end{array}
\right]
\]
It satisfies
\[
c\left(  X\right)  =x,c\left(  Y\right)  =y,c\left(  H\right)  =h.
\]

Now, as remarked earlier, the root system $\sum(\frak{t}_{\mathbb{C}}%
,\frak{g}_{\mathbb{C}})$ is of type $C_{n}$. Moreover, the various compact and
non-compact root multiplicities are as follows:
\begin{align*}
\dim\frak{k}_{\pm(\gamma_{i}\pm\gamma_{j})/2}  &  =d\text{, }\dim\frak{p}%
_{\pm(\gamma_{i}\pm\gamma_{j})/2}=d\\
\dim\frak{k}_{\pm\gamma_{i}}  &  =e\text{, }\dim\frak{p}_{\pm\gamma_{i}}=1
\end{align*}
We fix holomorphic homomorphisms
\[
\Psi_{j}:\frak{sl}_{2}(\mathbb{C})\longrightarrow\frak{g}_{\mathbb{C}}\text{,
}j=1,...,n
\]
such that each $\Psi_{j}(X)$ spans $\frak{p}_{\gamma_{j}}$, and we write
\[
X_{j}=\Psi_{j}(X),x_{j}=\Psi_{j}(x)\text{, }y_{j}=\Psi_{j}(y)\text{,etc.}%
\]
The images of $\Psi_{j}$ commute with each other and we also write
\[
\mathbf{X}=\sum_{j}X_{j}\text{, }\mathbf{x}=\sum_{j}x_{j},\mathbf{y}=\sum
_{j}y_{j},\text{etc.}%
\]
The \emph{Cayley transform} of $\frak{g}$ is the product
\[
\mathbf{c}=\exp\text{ ad }\frac{\pi i}{4}\left(  \mathbf{x}+\mathbf{y}\right)
=\exp\text{ ad }\frac{\pi i}{4}(\mathbf{X}+\mathbf{Y}).
\]
We write $\frak{a}=$ $\mathbf{c}\left(  i\frak{t}\right)  $ for the Cayley
transform of $i\frak{t}$. This is the abelian subalgebra of $\frak{g}$ spanned
by $h_{1},\cdots,h_{n}$.\label{sec-Striples}

\subsection{The symmetric space $L/H$}

Let $H\subset L$ be the stabilizer of $\mathbf{y}\in\overline{\frak{n}}$, then
condition 2$.$ implies (and is in fact equivalent to the assertion) that $L/H$
is a symmetric space. The involution \thinspace$\sigma$ for this symmetric
space consists of conjugation by a suitable element of $K$ --- corresponding
to the element $w_{0}$ of condition 2.

\textbf{Example. }If $G=O_{2n,2n}$, then $\,L=GL_{2n}(\mathbb{R})$ and $N$ is
the Jordan algebra of $2n\times2n$ real skew-symmetric matrices, and
$H=Sp_{n}(\mathbb{R})$.

In the present situation $L/H$ is always non-Riemannian; and moreover
$\frak{a}\ $is the corresponding Cartan subalgebra in the usual sense. In
other words, if we consider the Cartan decompositions for $\theta$ and
$\sigma$
\[
\frak{l}=\frak{m}+\frak{r}\text{ ,}\ \frak{l}=\frak{h}+\frak{q}\text{;}%
\]
then $\frak{a}$ is a Cartan subspace in $\frak{q\cap r}$.

Since $\frak{a}=\mathbf{c}\left(  i\frak{t}\right)  $, the roots of $\frak{a}$
in $\frak{g}$ are
\[
\Sigma(\frak{a},\frak{g})=\left\{  \pm\varepsilon_{i}\pm\varepsilon_{j}%
,\pm2\varepsilon_{j}\right\}  \text{ where }\varepsilon_{i}=\frac{1}{2}%
\gamma_{i}\ \circ\mathbf{c}^{-1}.
\]
Moreover it is easy to see that
\[
\Sigma(\frak{a},\frak{l})=\left\{  \pm(\varepsilon_{i}-\varepsilon
_{j})\right\}  ,\text{ }\Sigma(\frak{a},\frak{n})=\left\{  \varepsilon
_{i}+\varepsilon_{j},2\varepsilon_{j}\right\}  ,\text{ }\Sigma(\frak{a}%
,\overline{\frak{n}})=\left\{  -\varepsilon_{i}-\varepsilon_{j},-2\varepsilon
_{j}\right\}  .
\]

We now remark that the weight
\[
\nu=\varepsilon_{1}+\varepsilon_{2}+\ldots+\varepsilon_{n}%
\]
extends to a character of $\frak{l}$. The easiest way to see this is to
consider for $a$ in $\frak{a}$,
\[
\operatorname*{tr}\operatorname*{ad}\nolimits_{\overline{\frak{n}}}\left(
a\right)  =\left[  -2d\sum\left(  \varepsilon_{i}+\varepsilon_{j}\right)
-\left(  e+1\right)  \sum2\varepsilon_{j}\right]  \left(  a\right)
=-2r\nu\left(  a\right)  .
\]
Thus we can define
\[
\nu\left(  l\right)  =-\left(  \frac{1}{2r}\right)  \operatorname*{tr}%
\operatorname*{ad}\nolimits_{\overline{\frak{n}}}\left(  l\right)  \text{ for
}l\in\frak{l}\text{, where }r=d(n-1)+(e+1).\text{ }%
\]
Similarly, we define a corresponding positive character of $L$ by
\[
e^{\nu}\left(  l\right)  =l^{\nu}=\left|  \det\operatorname*{Ad}%
\nolimits_{\overline{\frak{n}}}l\right|  ^{-1/\left(  2r\right)  }\text{ for
}l\text{ in }L.
\]
Considering an appropriate power of the determinant we obtain a corresponding
positive character of the groups $P$ and $\overline{P}$, which we write as
$g\mapsto e^{\nu}\left(  g\right)  $, or as $g\mapsto g^{\nu}$.

To complete the connection with the Jordan structure, we note that the Jordan
norm $\phi$ on $\overline{\frak{n}}$ is a polynomial function which transforms
by the character $e^{-2\nu}$ of $L$. Finally we observe that the Killing form
on $\frak{g}$ gives a pairing between $\frak{n}$ and $\overline{\frak{n}}$
which we \emph{rescale} by setting $\left\langle x_{1},y_{1}\right\rangle =1$.\label{sec-LH}

\section{Integral formulas}

It is an immediate consequence of Corollary \ref{cor-equiv} that the
$L$-orbits in $\overline{\frak{n}}$ carry equivariant measures. More
precisely, we write $e^{\nu}$ for the positive character of $L$ defined in
\ref{sec-LH} and let $r=d(n-1)+(e+1)$ be as before. Then we have

\begin{lemma}
\label{lemma-orbitequiv}

\begin{enumerate}
\item  The Lebesgue measure $d\lambda$ on $\overline{\frak{n}}$ is $e^{2r\nu}
$-equivariant for the $L$-action.

\item  The rank $k$-orbit carries an $e^{2dk\nu}$-equivariant measure
$d\mu=d\mu_{k}$.
\end{enumerate}
\end{lemma}

The (easy) proof of this lemma is postponed to the next subsection. We now
describe a ``polar coordinates'' \ expression for these equivariant measures.
Let $\mathcal{O}$ be the rank $k$ orbit, $\mathcal{O}=L\cdot\left(
y_{1}+\ldots+y_{k}\right)  $. In \cite{loos} it is shown that the elements
\[
\{z_{1}y_{1}+\ldots+z_{k}y_{k}\,|\;z_{1}>z_{2}>\ldots>z_{k}>0\}
\]
give a complete set of orbit representatives for the action of $M=L\cap K$
\ on the rank $k$ orbit. Accordingly, we write $C_{k}\subset\mathbb{R}^{k}$
for the cone
\[
C_{k}=\{z=\left(  z_{1},z_{2},\ldots,z_{k}\right)  \,|\;z_{1}>z_{2}%
>\ldots>z_{k}>0\};
\]
and for $m$ in $M$, $z$ in $C_{k}$ we write
\[
m\cdot z=\operatorname*{Ad}m\left(  z_{1}y_{1}+\ldots+z_{k}y_{k}\right)
\in\overline{\frak{n}}.
\]

For $z$ in $C_{k}$ we introduce the notation
\[
P_{k}\left(  z\right)  =z_{1}\ldots z_{k},V_{k}\left(  z\right)  =\prod_{1\leq
i<j\leq k}\left[  z_{i}^{2}-z_{j}^{2}\right]  ,\text{ }d_{k}^{\times}%
z=\prod_{j=1}^{k}\frac{dz_{j}}{z_{j}}%
\]
where each $dz_{j}$ denotes the Lebesgue measure on $\mathbb{R}$.

The main results of this section are summarized in the following two propositions.

\begin{proposition}
\label{prop-open} Let $d\lambda$ be the Lebesgue measure on $\overline
{\frak{n}}$, then
\[
\int\limits_{\overline{\frak{n}}}fd\lambda=c\int\limits_{C_{n}}\left[
\int\limits_{M}f(m\cdot z)dm\right]  d_{\ast}z\text{\emph{\ where} }d_{\ast
}z=\left[  P_{n}\right]  ^{e+1}\left[  V_{n}\right]  ^{d}d_{n}^{\times}z.
\]
\end{proposition}

\begin{proposition}
\label{prop-sing} Let $d\mu$ be the equivariant measure on the rank $k$ orbit
$\mathcal{O}$, then
\[
\int\limits_{\mathcal{O}}fd\mu=c\int\limits_{C_{k}}\left[  \int\limits_{M}%
f(m\cdot z)dm\right]  d_{k}z\text{\emph{\ where} }d_{k}z=\left[  P_{k}%
^{n-k+1}V_{k}\right]  ^{d}d_{k}^{\times}z.
\]
\end{proposition}

The scalars $c$ appearing in the above integral formulas are independent of
$f$ and depend only on the normalization of the measures $d\lambda$ and $d\mu$.

For subsequent purposes we need to consider the Lebesgue measure on $\frak{n}$
as well as $\overline{\frak{n}}$. For $m$ in $M$, $z$ in $C_{n},$ we write
\begin{equation}
m\circ z=\operatorname*{Ad}m\left(  z_{1}x_{1}+\ldots+z_{n}x_{n}\right)
\in\frak{n\,}. \label{=mz}%
\end{equation}
Since the Cartan involution $\theta$ gives a linear isomorphism between
$\frak{n}$ and $\overline{\frak{n}}$, satisfying $\theta\left(  m\circ
z\right)  =m\cdot z$, the following result can be derived immediately from
Proposition \ref{prop-open}.

\begin{corollary}
\label{cor-open} Let $d\lambda$ be the Lebesgue measure on $\frak{n}$, then
for all functions $f$ on $\frak{n}$
\[
\int\limits_{\frak{n}}fd\lambda=c\int\limits_{C_{n}}\left[  \int
\limits_{M}f(m\circ z)dm\right]  d_{\ast}z\text{\emph{\ where} }d_{\ast
}z=\left[  P_{n}\right]  ^{e+1}\left[  V_{n}\right]  ^{d}d_{n}^{\times}z.
\]
\endproof
\end{corollary}

The proofs of the propositions will occupy the rest of this section.

\subsection{Stabilizers and equivariant measures}

In this subsection, we prove Lemma \ref{lemma-orbitequiv}. For this we need to
first determine the stabilizer of the point $y_{1}+\ldots+y_{k}$ in the rank
$k$ orbit. If $k=n$, the stabilizer is the symmetric subgroup $H$ described
previously. We now discuss the remaining orbits. To simplify notation we
fix$\ k$ and write
\[
\mathbf{y}^{1}=y_{1}+\ldots+y_{k}.
\]
In Jordan algebra terms, $\mathbf{y}^{1}$ is a Peirce idempotent and
considering the $1$ and $0$ Peirce-eigenspaces of $\mathbf{y}^{1}$, we obtain
smaller Jordan algebras $\overline{\frak{n}}_{1}$ and $\overline{\frak{n}}%
_{0}$ with identity elements $\mathbf{y}^{1}$ and $\mathbf{y}^{0}%
=y_{k+1}+\ldots+y_{n}$, respectively. The corresponding structure groups
$L_{1}$ and $L_{0}$ are naturally subgroups of $L$. Subgroups of $L_{1}$ and
$L_{0}$ will be distinguished by subscripts $1$ and $0$, respectively. For
example,
\[
M_{1}=M\cap L_{1}\text{, }M_{0}=M\cap L_{0}\text{, }H_{1}=H\cap L_{1}\text{,
}\frak{a}_{1}=\frak{a}\cap\frak{l}_{1}.
\]
Thus $H_{1}$ is the stabilizer of $\mathbf{y}^{1}$ in $L_{1}$, and the full
stabilizer of $\mathbf{y}^{1}$ in $L$ is given by
\begin{equation}
S=(H_{1}\times L_{0})\cdot\overline{U}, \label{=stabilizer}%
\end{equation}
where $\overline{U}$ is the abelian subgroup whose Lie algebra $\overline
{\frak{u}}$ is spanned by the root spaces $\frak{l}^{-\varepsilon
_{i}+\varepsilon_{j}}$ ($1\leq i\leq k<j\leq n$). Then
\[
\frak{l}=\frak{s+}\left(  \frak{q}_{1}+\frak{u}\right)  \text{, where
}\frak{q}_{1}=\frak{q}\cap\frak{l}_{1},\frak{u}=\theta\overline{\frak{u}}.
\]

\begin{proof}
(of Lemma \ref{lemma-orbitequiv}) From the calculation in the previous
section, it follows that the Lebesgue measure $d\lambda$ on $\overline
{\frak{n}}$ is equivariant by the character
\[
\eta\left(  l\right)  =\left|  \det\operatorname*{Ad}\nolimits_{\overline
{\frak{n}}}l\right|  ^{-1}=e^{2r\nu}\left(  l\right)  =l^{2\left[
d(n-1)+e+1\right]  \nu}\text{ for }l\text{ in }L\text{,}%
\]
which proves the first part of the lemma.

For the second part, we consider the stabilizer of $y_{1}+\ldots+y_{k}$. First
suppose $k=n$; in this case the stabilizer $H$ is semisimple, and hence
\[
\left|  \det\operatorname*{Ad}\nolimits_{\frak{h}}h\right|  =1\text{ for all
}h\in H.
\]
Since $e^{\nu}|_{H}=1$ as well, it follows from Corollary \ref{cor-equiv} that
for each real $t$, the open orbit carries a measure which is $e^{t\nu}%
$-equivariant. The measure $d\mu$ is simply the special case $t=2dn$.

Now suppose $k<n$, then $\frak{s}$ $=\frak{h}_{1}+\frak{l}_{0}+\overline
{\frak{u}}$ and it is easy to see that $\operatorname*{tr}\left[
\operatorname*{ad}\nolimits_{\frak{s}}\left(  z\right)  \right]  =0$ for $z$
in $\frak{h}_{1}$ and $\frak{u}$. Thus we can assume that $z$ is in
$\frak{l}_{0}$, and then $\operatorname*{tr}\left[  \operatorname*{ad}%
\nolimits_{\frak{h}_{1}+\frak{l}_{0}}\left(  z\right)  \right]  =0,$ hence
$\operatorname*{tr}\left[  \operatorname*{ad}\nolimits_{\frak{s}}\left(
z\right)  \right]  =\operatorname*{tr}\left[  \operatorname*{ad}%
\nolimits_{\overline{\frak{u}}}\left(  z\right)  \right]  $. To calculate
this, it suffices to consider $z$ in $\frak{a}_{0}$, and then we get
\begin{align*}
\operatorname*{tr}\left[  \operatorname*{ad}\nolimits_{\frak{s}}\left(
z\right)  \right]   &  =\operatorname*{tr}\left[  \operatorname*{ad}%
\nolimits_{\overline{\frak{u}}}\left(  z\right)  \right]  =\sum_{1\leq i\leq
k<j\leq n}2d\left(  -\varepsilon_{i}+\varepsilon_{j}\right)  \left(  z\right)
\\
&  =-2d\left(  n-k\right)  \left[  \varepsilon_{1}+\ldots+\varepsilon
_{k}\right]  \left(  z\right)  +2dk\left[  \varepsilon_{k+1}+\ldots
+\varepsilon_{n}\right]  \left(  z\right)  .
\end{align*}
Now the weights $\varepsilon_{1},\ldots,\varepsilon_{k}$ restrict trivially on
$\frak{a}_{0}$, hence we see that the required trace is the restriction of
$\ 2dk\nu$ to $\frak{s}$. Thus we get
\[
\left|  \det\operatorname*{Ad}\nolimits_{\frak{s}}s\right|  =e^{2dk\nu}\left(
s\right)
\]
and the result follows by Corollary \ref{cor-equiv}.
\end{proof}

\label{sec-equiv}

\subsection{Jacobians for homogeneous spaces}

In order to prove integral formulas for homogeneous spaces, we need a method
for calculate the Jacobian of a diffeomorphism between such spaces. It is
convenient to work with the more flexible notion of a \emph{local
diffeomorphisms} between $X$ and $Y$, by which we simply mean a diffeomorphism
between open subsets of $X$ and $Y$. We also introduce the notation
\[
\left(  X,x\right)  \overset{F}{\rightarrow}\left(  Y,y\right)
\]
to represent the situation where $X$ and $Y$ are smooth manifolds; $x\ $\ and
$y$ are points in $X$ and $Y$; $F$ is a diffeomorphism from an open
neighborhood of $x$ to an open neighborhood of $Y$ such that $F\left(
x\right)  =y$.

We now consider the following situation: Suppose $X$ and $Y$ are homogenous
spaces for groups $G$ and $H$, and $dx$ and $dy$ are regular measures on this
spaces which are equivariant for characters $\gamma\left(  g\right)  \ $and
$\eta\left(  h\right)  $ respectively. We choose two points $x_{0}\in X$ and
$y_{0}\in Y$, and fix linear bases for the tangent spaces at these points.

\begin{lemma}
\label{lemma-jaccalc} Suppose that in the above situation we have a local
diffeormorphism $F:X\rightarrow Y$. Then for all $x$ in the domain of $F$ we
have
\[
J_{F}\left(  x\right)  =c\eta\left(  h\right)  \gamma\left(  g\right)  \left|
\det D_{hFg}\left(  x_{0}\right)  \right|  ,
\]
where $c$ is a scalar independent of $x$; $h\in H$ and $g\in G$ satisfy
\[
g\cdot x_{0}=x\text{, }h\cdot F\left(  x\right)  =y_{0}%
\]
and $D_{hFg}:T_{x_{0}}X\rightarrow T_{y_{0}}Y$ is regarded as a matrix via the
above bases.
\end{lemma}

\begin{proof}
We fix $x$ and write $y=F\left(  x\right)  $. Then we have local
diffeomorphisms
\[
\left(  X,x_{0}\right)  \overset{g}{\rightarrow}\left(  X,x\right)
\overset{F}{\rightarrow}\left(  Y,y\right)  \overset{h}{\rightarrow}\left(
Y,y_{0}\right)  .
\]
By formula (\ref{prodjac}) we get
\[
J_{hFg}\left(  x_{0}\right)  =J_{g}\left(  x_{0}\right)  J_{F}\left(
x\right)  J_{h}\left(  y\right)  =\gamma\left(  g\right)  ^{-1}J_{F}\left(
x\right)  \eta\left(  h\right)  ^{-1}.
\]
On the other hand by Lemma \ref{lemma-jacdet} we have
\[
J_{hFg}\left(  x_{0}\right)  =c\left|  \det D_{hFg}\left(  x_{0}\right)
\right|
\]
for some positive scalar $c$, independent of $hFg$. The lemma follows.
\end{proof}

\subsection{Integral formula for the Lebesgue measure}

We now apply the results of the previous subsections to prove the integration
formulas on the $L$-orbits. In order to do this we first need to fix bases for
various subspaces of $\frak{l}$ which are compatible with the actions of
$\theta$ and $\sigma$. We start with the weight decomposition
\[
\frak{l}=\frak{l}^{0}\oplus\left(  \oplus_{i\neq j}\frak{l}^{\varepsilon
_{i}-\varepsilon_{j}}\right)
\]
where $\frak{l}^{0}$ is the centralizer of $\frak{a}$. The involutions
$\theta$ and $\sigma$ act by $-1$ on $\frak{a}$, thus the space $\frak{l}^{0}$
is invariant by $\theta$ and $\sigma$; considering their eigenvalues we have
the decomposition
\[
\frak{l}^{0}=\frak{l}^{++}+\frak{l}^{+-}+\frak{l}^{-+}+\frak{l}^{--}%
\]
where $\frak{l}^{++}=\frak{l}^{0}\cap(\frak{m}\cap\frak{h}),\frak{l}%
^{+-}=\frak{l}^{0}\cap(\frak{m}\cap\frak{q})$ etc. Note of course that
$\frak{l}^{--}=\frak{a}$. On the other hand, the root spaces $\frak{l}%
^{\alpha}$ are not stable under $\theta$ and $\sigma$, in fact each involution
maps $\frak{l}^{\alpha}$ to $\frak{l}^{-\alpha}$. However the involution
$\tau=\sigma\theta=\theta\sigma$ does stabilize these spaces. Considering the
eigenvalues of $\tau$ we have a decomposition
\[
\frak{l}^{\alpha}=\frak{l}^{\alpha,+}+\frak{l}^{\alpha,-}%
\]
where $\frak{l}^{\alpha,+}=\frak{l}\cap(\frak{m}\cap\frak{h}+\frak{r}%
\cap\frak{q})$, $\frak{l}^{\alpha,-}=\frak{l}^{\alpha}\cap(\frak{m}%
\cap\frak{q}+\frak{r}\cap\frak{h})$ (cf. \cite[8.1]{sch}).

In the present situation, we \emph{always\ }have
\[
\dim\frak{l}^{\alpha,+}=\text{ }\dim\frak{l}^{\alpha,-}=d\text{.}%
\]
This fact can be checked easily for each of the examples in the table below,
from the lists of multiplicities in \cite{oshima}. We fix bases
\[
\left\{  X_{l}^{\alpha,\pm}\right\}  \subset\frak{l}^{\alpha,\pm}\text{,
}\left\{  X_{l}^{\pm,\pm}\right\}  \subset\frak{l}^{\pm,\pm}\text{.}%
\]

We now turn to the proof of Proposition \ref{prop-open}. The rank $n$ orbit
$L/H$ is open and dense in $\overline{\frak{n}}$ and its complement has
measure $0$. Since $L/H$ is a symmetric space of type $A_{n-1}$, the
``multiplication'' map
\[
F:\overline{m}\times a\mapsto\overline{ma^{-1}}%
\]
gives an diffeomorphism between $\left[  M/M^{\prime}\right]  \times A^{+}$
and $L/H$, where $M^{\prime}$ is the centralizer of $A$ in $M\cap H$ and
$A^{+}=\exp$ $\frak{a}^{+}$ with
\[
\frak{a}^{+}=\left\{  c_{1}h_{1}+\cdots+c_{n}h_{n}\mid c_{1}\geq\ldots\geq
c_{n}\right\}
\]
We regard $F$ as a local diffeomorphism between the homogeneous spaces
$\left[  M/M^{\prime}\right]  \times A^{+}=\left[  M\times A\right]  /\left[
M^{\prime}\times\left\{  1\right\}  \right]  $ and $L/H$. By Corollary
\ref{cor-equiv}, the first space carries an invariant measure $d\overline
{m}\times da$, while the second space carries the Lebesgue measure $d\lambda$
which is $e^{2r\nu}$-equivariant by Lemma \ref{lemma-orbitequiv}. The main
result is the following Jacobian computation.

\begin{lemma}
\label{lemma-openjac} In the above situation there is a scalar $c$ such that
\[
J_{F,d\overline{m}\times da,d\lambda}(\overline{m}\times a)=ca^{2\left(
e+1\right)  \nu}\prod_{i<j}\left|  a^{4\varepsilon_{i}}-a^{4\varepsilon_{j}%
}\right|  ^{d}.
\]
\end{lemma}

\begin{proof}
We fix the natural base points $x_{0}=\overline{1}\times1\in$ $\left[
M/M^{\prime}\right]  \times A$ and $y_{0}=\overline{1}\in L/H$ and apply Lemma
\ref{lemma-jaccalc} from the previous subsection with
\[
g=m\times a,h=am^{-1}%
\]
Then we have
\[
hFg\left(  \overline{m^{\prime}},a^{\prime}\right)  =hF\left(  \overline
{mm^{\prime}},aa^{\prime}\right)  =h\left(  \overline{mm^{\prime}\left(
aa^{\prime}\right)  ^{-1}}\right)  =\overline{\left(  am^{\prime}%
a^{-1}\right)  a^{\prime-1}}%
\]
Therefore the differential $D_{hFg}:\frak{m}/\frak{m}^{\prime}%
+\frak{a\rightarrow l}/\frak{h}$ is given by
\begin{equation}
T_{a}\left(  v,w\right)  :=\operatorname*{Ad}\left(  a\right)  v+w\text{
}\left(  \operatorname{mod}\text{ }\frak{h}\right)  \label{=Ta}%
\end{equation}
This shows that
\[
J_{F,d\overline{m}\times da,d\lambda}(\overline{m}\times a)\ \sim a^{2r\nu
}\left|  \det\text{ }T_{a}\right|
\]
for some fixed choice of bases for $\frak{m}/\frak{m}^{\prime}+\frak{a\ }$and
$\frak{l}/\frak{h}$. \ Since $r=d(n-1)+e+1$, the result follows from the next lemma.
\end{proof}

\begin{lemma}
\label{lemma-Ta} For the map $T_{a}$ defined in formula (\ref{=Ta}) we have
\[
\left|  \det\text{ }T_{a}\right|  \sim a^{-2d\left(  n-1\right)  \nu}%
\prod_{i<j}\left|  a^{4\varepsilon_{i}}-a^{4\varepsilon_{j}}\right|  ^{d}%
\]
\end{lemma}

\begin{proof}
In order to compute the determinant, we choose convenient bases for
$\frak{l}/\frak{h}$ and $\frak{m}/\frak{m}^{\prime}+\frak{a}$. We have
$\frak{l}/\frak{h}\approx\frak{q}$, and so\ we may choose the basis
\[
\left\{  X_{k}^{\pm,-}\right\}  ,\left\{  \left[  X_{k}^{\alpha,+}-\theta
X_{k}^{\alpha,+}\right]  ,\left[  X_{k}^{\alpha,-}+\theta X_{k}^{\alpha
,-}\right]  \mid\alpha>0\right\}  ;
\]

For the space $\frak{m}/\frak{m}^{\prime}+\frak{a}$, we note that
$\frak{m}^{\prime}=\frak{l}^{++}$ in our earlier notation, and thus a basis is
given by:
\[
\left\{  X_{k}^{\pm,-}\right\}  ,\left\{  \left[  X_{k}^{\alpha,+}+\theta
X_{k}^{\alpha,+}\right]  ,\left[  X_{k}^{\alpha,-}+\theta X_{k}^{\alpha
,-}\right]  \mid\alpha>0\right\}
\]

We now claim that
\begin{align*}
T_{a}\cdot X_{k}^{\pm,-}  &  =X_{k}^{\pm,-}\\
T_{a}\cdot\left[  X_{k}^{\alpha,+}+\theta X_{k}^{\alpha,+}\right]   &
=\frac{a^{\alpha}-a^{-\alpha}}{2}\left[  X_{k}^{\alpha,+}-\theta X_{k}%
^{\alpha,+}\right]  \text{ }\left(  \operatorname{mod}\text{ }\frak{h}\right)
\\
T_{a}\cdot\left[  X_{k}^{\alpha,-}+\theta X_{k}^{\alpha,-}\right]   &
=\frac{a^{\alpha}+a^{-\alpha}}{2}\left[  X_{k}^{\alpha,-}+\theta X_{k}%
^{\alpha,-}\right]  \text{ }\left(  \operatorname{mod}\text{ }\frak{h}\right)
\end{align*}
The first equality is obvious. For the second, we calculate
\begin{align*}
T_{a}\cdot\left[  X_{k}^{\alpha,+}+\theta X_{k}^{\alpha,+}\right]   &
=\operatorname*{Ad}\left(  a\right)  \cdot\left[  X_{k}^{\alpha,+}+\theta
X_{k}^{\alpha,+}\right]  =a^{\alpha}X_{k}^{\alpha,+}+a^{-\alpha}\theta
X_{k}^{\alpha,+}\\
&  =\frac{a^{\alpha}-a^{-\alpha}}{2}\left[  X_{k}^{\alpha,+}-\theta
X_{k}^{\alpha,+}\right]  +\frac{a^{\alpha}+a^{-\alpha}}{2}\left[
X_{k}^{\alpha,+}+\theta X_{k}^{\alpha,+}\right] \\
&  =\frac{a^{\alpha}-a^{-\alpha}}{2}\left[  X_{k}^{\alpha,+}-\theta
X_{k}^{\alpha,+}\right]  \text{ }\left(  \operatorname{mod}\text{ }%
\frak{h}\right)  .
\end{align*}
The third equality follows by a similar calculation.

Since there are $d$ vectors in each of the sets $\left\{  X_{k}^{\alpha
,+}+\theta X_{k}^{\alpha,+}\right\}  ,\left\{  X_{k}^{\alpha,-}+\theta
X_{k}^{\alpha,-}\right\}  $ we get
\begin{align*}
\left|  \det T_{a}\right|   &  \sim\prod_{i<j}\left|  \frac{a^{\varepsilon
_{i}-\varepsilon_{j}}-a^{\varepsilon_{j}-\varepsilon_{i}}}{2}\right|
^{d}\left|  \frac{a^{\varepsilon_{i}-\varepsilon_{j}}+a^{\varepsilon
_{j}-\varepsilon_{i}}}{2}\right|  ^{d}\\
&  \sim\prod_{i<j}\left|  a^{2\varepsilon_{i}-2\varepsilon_{j}}%
-a^{2\varepsilon_{j}-2\varepsilon_{i}}\right|  ^{d}=\prod_{i<j}\left|
a^{-2\left(  \varepsilon_{i}+\varepsilon_{j}\right)  }\left(  a^{4\varepsilon
_{i}}-a^{4\varepsilon_{j}}\right)  \right|  ^{d}\\
&  =a^{-2d\left(  n-1\right)  \nu}\prod_{i<j}\left|  a^{4\varepsilon_{i}%
}-a^{4\varepsilon_{j}}\right|  ^{d}.
\end{align*}
\end{proof}

We can now prove Proposition \ref{prop-open}.

\begin{proof}
(of Proposition \ref{prop-open}) Since the map $F$ is a diffeomorphism between
$\left[  M/M^{\prime}\right]  \times A^{+}$ and $L/H,$ by the Jacobian
calculation, we get
\[
\int_{\overline{\frak{n}}}fd\lambda=\int_{L/H}fd\lambda=c\int_{A^{-}}\left[
\int_{M/M^{\prime}}f\left(  \overline{ma^{-1}}\right)  d\overline{m}\right]
a^{2\left(  e+1\right)  \nu}\prod_{i<j}\left|  a^{4\varepsilon_{i}%
}-a^{4\varepsilon_{j}}\right|  ^{d}da.
\]
For a suitable normalization of the Haar measures, the innermost integral can
be rewritten as $\int_{M}f\left(  \overline{ma^{-1}}\right)  dm$. Also for
$a=\exp\left(  c_{1}h_{1}+\cdots+c_{n}h_{n}\right)  $ in $A^{+}$, we have
\[
\overline{a^{-1}}=a\cdot\mathbf{y=}e^{2c_{1}}y_{1}+\cdots+e^{2c_{n}}y_{n}.
\]
Thus if we identify $A^{+}$ and $C_{n}$ via the map
\[
\exp\left(  c_{1}h_{1}+\cdots+c_{n}h_{n}\right)  \longleftrightarrow\left(
e^{2c_{1}},\cdots,e^{2c_{n}}\right)
\]
then we get
\[
\overline{ma^{-1}}\leftrightarrow m\cdot z,\,\,a^{2\nu}\leftrightarrow
P_{n},\,\,\prod_{i<j}\left|  a^{4\varepsilon_{i}}-a^{4\varepsilon_{j}}\right|
\leftrightarrow V_{n},\,\,da\leftrightarrow d_{n}^{\times}z
\]
Making these substitutions in the integral formula, we obtain the desired result.
\end{proof}

\subsection{Integral formula for the equivariant measure}

We now turn to the proof of Proposition \ref{prop-sing}. For the case of the
measure $d\mu$ on the open orbit, the calculations are \emph{exactly }the same
as in the previous section, the only difference being that the measure is
$e^{2dn\nu}$-equivariant. Thus arguing as in Lemmas \ref{lemma-openjac} and
\ref{lemma-Ta} we obtain that
\[
J_{F,d\overline{m}\times da,d\mu}(\overline{m}\times a)\sim a^{2dn\nu}\left|
\det T_{a}\right|  =a^{2d\nu}\prod_{i<j}\left|  a^{4\varepsilon_{i}%
}-a^{4\varepsilon_{j}}\right|  ^{d}%
\]
Now arguing as in the proof of Proposition \ref{prop-open} we obtain
Proposition \ref{prop-sing} for $k=n$.

Now consider the case $k\leq n$. We recall the notation $\mathbf{y}^{1}%
=y_{1}+\ldots+y_{k},M_{1},M_{0}$, etc. introduced in section \ref{sec-equiv}.
Now writing $M_{1}^{\prime}=Z_{M_{1}\cap H_{1}}(\frak{a}_{1})$ and
$\frak{a}_{1}^{+}=\frak{a}_{1}\cap\frak{a}^{+}$ etc., we get a local
diffeomorphism $F^{1}:M/(M_{1}^{\prime}\times M_{0})\times A_{1}%
^{+}\rightarrow\mathcal{O}=L/S$ given by
\[
F^{1}(\overline{m},a)=\overline{ma^{-1}}%
\]
Arguing as in Lemma \ref{lemma-openjac}, we obtain
\[
J_{F,d\overline{m}\times da,d\mu}(\overline{m}\times a)\sim a^{2dk\left(
\varepsilon_{1}+\cdots+\varepsilon_{k}\right)  }\left|  \det T_{a}%
^{1}\right|
\]
where $T_{a}^{1}:\frak{m}/\left(  \frak{m}_{1}^{\prime}+\frak{m}_{0}\right)
+\frak{a}_{1}\rightarrow\frak{l}/\frak{s}$ is given by
\[
T_{a}^{1}\left(  v,w\right)  =\operatorname*{Ad}a\left(  v\right)  +w\text{
}\left(  \operatorname{mod}\frak{s}\text{ }\right)
\]
and we compute the determinant for some choice of basis for the two sides. We
choose these bases in the following manner. We have $A_{1}$-module
isomorphisms
\begin{align*}
\frac{\frak{m}}{\frak{m}_{1}^{\prime}+\frak{m}_{0}}+\frak{a}_{1}  &
\approx\frac{\frak{m}}{\frak{m}_{1}+\frak{m}_{0}}+\frac{\frak{m}_{1}}%
{\frak{m}_{1}^{\prime}}+\frak{a}_{1}\text{ }\\
\frak{l}/\frak{s}  &  \approx\frak{u}+\frak{q}_{1}.
\end{align*}
Thus the map $T_{a}^{1}$ decomposes as a direct sum $T_{a}^{1}=T_{0}+T_{1}$
with
\[
T_{0}:\frac{\frak{m}}{\frak{m}_{1}+\frak{m}_{0}}\rightarrow\frak{u}\text{,
}T_{1}:\frac{\frak{m}_{1}}{\frak{m}_{1}^{\prime}}+\frak{a}_{1}\rightarrow
\frak{q}_{1}.
\]
The calculation for $\left|  \det T_{1}\right|  $ is the same as in Lemma
\ref{lemma-Ta}, applied to a smaller Jordan algebra of rank $k$. This gives
\[
\left|  \det T_{1}\right|  =\prod_{i<j\leq k}\left|  a^{2\varepsilon
_{i}-2\varepsilon_{j}}-a^{2\varepsilon_{j}-2\varepsilon_{i}}\right|
^{d}=a^{-2d\left(  k-1\right)  \left(  \varepsilon_{1}+\cdots+\varepsilon
_{k}\right)  }\prod_{i<j\leq k}\left|  a^{4\varepsilon_{i}}-a^{4\varepsilon
_{j}}\right|  ^{d}%
\]
To calculate $\det T_{0}$ we use the following bases for $\frak{m}/\left(
\frak{m}_{1}+\frak{m}_{0}\right)  $ and $\frak{u}$
\begin{align*}
\frac{\frak{m}}{\frak{m}_{1}+\frak{m}_{0}}  &  :\left\{  \left[  X_{l}%
^{\alpha,\pm}+\theta X_{l}^{\alpha,\pm}\right]  \mid\alpha=\varepsilon
_{i}-\varepsilon_{j},\text{ }i\leq k,j>k\right\}  .\\
\frak{u}  &  :\left\{  X_{l}^{\alpha,\pm}\mid\alpha=\varepsilon_{i}%
-\varepsilon_{j},\text{ }i\leq k,j>k\right\}
\end{align*}
The action of $T_{0}$ is given by
\[
\operatorname*{ad}a\cdot\left[  X_{l}^{\alpha,\pm}+\theta X_{l}^{\alpha,\pm
}\right]  =a^{\alpha}X_{l}^{\alpha,\pm}+a^{-\alpha}\theta X_{l}^{\alpha,\pm
}=a^{\alpha}X_{l}^{\alpha,\pm}\text{ }\left(  \operatorname{mod}\text{
}\frak{s}\right)
\]
Since for $a$ in $A_{1}$ we have $a^{\varepsilon_{j}}=1$ for $j>k$, we get
\[
\left|  \det T_{0}\right|  =\prod_{i\leq k,j>k}\left|  a^{\varepsilon
_{i}-\varepsilon_{j}}\right|  ^{2d}=a^{2d\left(  n-k\right)  \left(
\varepsilon_{1}+\cdots+\varepsilon_{k}\right)  }.
\]
Thus we get
\begin{align*}
J_{F,d\overline{m}\times da,d\mu}(\overline{m}\times a)  &  \sim a^{2dk\left(
\varepsilon_{1}+\cdots+\varepsilon_{k}\right)  }\left|  \det T_{1}\right|
\left|  \det T_{0}\right| \\
&  =a^{2d\left(  n-k+1\right)  \left(  \varepsilon_{1}+\cdots+\varepsilon
_{k}\right)  }\prod_{i<j\leq k}\left|  a^{4\varepsilon_{i}}-a^{4\varepsilon
_{j}}\right|  ^{d}.
\end{align*}

Now arguing as in the proof of Proposition \ref{prop-open} we obtain
Proposition \ref{prop-sing} for $k<n$.

\section{Estimates for spherical vectors}

We can relate the $P$-representation $\pi_{\mathcal{O}}$ to a unitarizable
submodule of a certain degenerate principal series for $G$, which is described
as follows: If $\chi$ is a character of $L$, we write $\left(  \pi_{\chi
},I(\chi)\right)  $ for the degenerate principal series representation
$\operatorname{Ind}{}_{\overline{P}}^{G}\chi$ (unnormalized smooth induction);
thus
\[
I(\chi)=\{f\in C^{\infty}(G)|f(l\overline{n}g)=\chi(l)f(g)\text{ for }l\in
L,\overline{n}\in\overline{N},g\in G\}
\]
and the group $G$ acts by right translations. By virtue of the Gelfand-Naimark
decomposition $G\approx\overline{P}N$ , functions from $I(\chi)$ are
determined by their restriction to $N$. Combining this with the exponential
map we can identify $I(\chi)$ with a subspace $E(\chi)$ of $C^{\infty
}(\frak{n})$\thinspace. We refer to this as the noncompact picture.

For $t\in\mathbb{R}$, we write $I(t)$, $E(t)$ for $I(e^{t\nu})$, $E(e^{t\nu})
$; more generally, if $\varepsilon:L\rightarrow\mathbb{T}$ is a unitary
character, we write $I(t,\varepsilon)$, $E(t,\varepsilon)$ for $I(e^{t\nu
}\otimes\varepsilon)$, $E(e^{t\nu}\otimes\varepsilon)$. These principal series
were studied in \cite{sahi-dp} via the ``Cayley operator'' $D$ which is the
constant coefficient differential operator on $\frak{n}$, whose symbol is the
Jordan norm polynomial $\phi$. Powers of $D$ are intertwining operators for
the principal series, and their eigenvalues on the various $K$-isotypic
components are given by the Capelli identity of \cite{kostant-sahi}.

$E(t)$ is a spherical representation of $G$ and we write $\Phi_{t}$ for the
$K$-spherical vector. Among the results obtained in \cite{sahi-dp} is that
for
\[
\chi_{k}=e^{-dk\nu},\text{ }k=1,\ldots,n-1;
\]
the space $E(-dk)$ contains a unitarizable spherical submodule. We need to
study the Fourier transforms of the corresponding spherical vectors
\[
\Phi_{-dk};\text{ }k=1,\ldots,n-1.
\]
For this we identify $\overline{\frak{n}}$ with the dual of $\frak{n}^{\ast}$
via the normalized Killing form from subsection \ref{sec-LH}. Also we fix
$k<n$, write $\Phi$ for the spherical vector $\Phi_{-dk}$, and write $\left(
\mathcal{O},d\mu\right)  $ for the rank $k$ orbit in $\overline{\frak{n}}$
together with its equivariant measure described in Lemma
\ref{lemma-orbitequiv}. The main results are

\begin{proposition}
\label{prop-L1}The measure $\Phi d\lambda$ is a tempered distribution on
$\frak{n}$ and there exists an $M$-invariant function $g$ in $\mathcal{L}%
^{1}(\mathcal{O},d\mu)$ such that
\[
\Phi d\lambda=\widehat{gd\mu}\text{. }%
\]
\end{proposition}

\begin{proposition}
For $k<n$, one has $g\in\mathcal{L}^{2}(\mathcal{O},d\mu)$. \label{=prop-L2}
\end{proposition}

We prove these propositions in the next few subsections. The strategy is as
follows: Let us write $\Phi_{k,n}$ for the function $\Phi_{-dk}$, in order to
emphasize dependence on $n$ as well as $k$. Now although the above results are
\emph{false in general }for the open orbit ($k=n$), nevertheless we can prove
the desired results by reducing to a slightly weaker estimate for $k=n$, which
turns out to be true, and somewhat easier to prove. We establish this result
in the next subsection and then outline the reduction procedure in the two
following subsections. We then deduce Theorem \ref{theoremA} from Proposition
\ref{=prop-L2} by arguments similar to \cite{sahi-expl} and \cite{hilbert-two}.\label{sec-dp}

\subsection{Estimates for the open orbit}

As indicated above, we first consider the function
\[
\Phi=\Phi_{n,n}=\Phi_{-dn}.
\]
We need appropriate $\mathcal{L}^{2}$-estimates with respect to the Lebesgue
measure $d\lambda$ on $\frak{n}$ for the function $\Psi$ and its derivatives.
The ``straightforward'' estimate is actually false for the group
$Sp_{n}(\mathbb{C})$, but it does work for the other groups $G$ in the table
in subsection \ref{sec-table}. Thus we formulate two results, one for $G$
$\neq Sp_{n}(\mathbb{C})$ and the other for all groups:

\begin{proposition}
For all groups $G$ \emph{other than} $Sp_{n}(\mathbb{C})$, we have
\thinspace$\Phi\in\mathcal{L}^{2}(\frak{n},d\lambda)$. \label{=lemma-phi0}
\end{proposition}

\begin{proposition}
For all groups $G$ and for all $m\geq1$, we have \thinspace$D^{m}\Phi
\in\mathcal{L}^{2}(\frak{n},d\lambda)$. \label{=lemma-dmphi}
\end{proposition}

For each $t$, the function $\Phi_{t}$ is $M$-invariant, and is therefore
determined by the restriction to the subspace $\left\{  z_{1}x_{1}%
+\cdots+z_{n}x_{n}\right\}  \subseteq\frak{n}$; we start by giving an explicit
formula for the restriction.

\begin{lemma}
We have
\[
\Phi_{t}(m\circ z)=\prod_{i=1}^{n}(1+z_{i}^{2})^{\frac{t}{2}}\text{ for all
}m\text{ in }M.
\]
\label{=lemma-growth}
\end{lemma}

\begin{proof}
For the group $G=SL_{2}\left(  \mathbb{R}\right)  $ this is a straightforward
calculation which we leave to the reader. In the general case, we view $\Phi$
as a function on $G$ which is right $K$-invariant, and left $\overline{P}%
$-equivariant with character $e^{t\nu}$. We now restrict $\Phi$ to the
subgroup $SL_{2}\times\cdots\times SL_{2}$ corresponding to the $S$-triples of
section \ref{sec-Striples}. This restriction is right $SO_{2}\times
\cdots\times SO_{2}$-invariant, and left $\overline{B}\times\cdots
\times\overline{B}$-equivariant with character $e^{s\nu}=$ $e^{s\varepsilon
_{1}}\times\cdots\times e^{s\varepsilon_{n}}$. Thus applying the $SL_{2}%
$-calculation to each factor, we conclude that the restriction to $z_{1}%
x_{1}+\cdots+z_{n}x_{n}$ is given as in the statement of the lemma .
\end{proof}

Combining this with Corollary \ref{cor-open} we obtain the following estimate

\begin{lemma}
For $t<-\left[  d\left(  n-1\right)  +\left(  e+1\right)  /2\right]  $, we
have $\Phi_{t}$ $\in\mathcal{L}^{2}(\frak{n},d\lambda)$. \label{lemma-intcond}
\end{lemma}

\begin{proof}
Combining the previous lemma with Corollary \ref{cor-open}, we get
\[
\int\left|  \Phi_{t}\right|  ^{2}d\lambda=\int_{C_{n}}\psi(z)dz_{1}%
dz_{2}\ldots dz_{n},
\]
where
\[
\psi(z)=\prod_{i=1}^{n}z_{i}^{e}(1+z_{i}^{2})^{t}\prod_{1\leq i<j\leq n}%
(z_{i}^{2}-z_{j}^{2})^{d}.
\]
Expanding $(z_{i}^{2}-z_{j}^{2})^{d}$, we can write $\psi(z)$ as a combination
of terms of the form
\[
\prod_{i=1}^{n}z_{i}^{e+k_{i}}(1+z_{i}^{2})^{t}\text{ where each }k_{i}%
\leq2d\left(  n-1\right)  .
\]
Each of these integrals is a \emph{product of one-variable} integrals which
converge if
\[
\int_{0}^{\infty}x^{e+2d\left(  n-1\right)  }\left(  1+x^{2}\right)
^{t}dx<\infty.
\]
This happens if $\ 2t+e+2d\left(  n-1\right)  <-1$, which proves the lemma.\smallskip
\end{proof}

\begin{corollary}
If $f\in E(t,\varepsilon)$ for some $t<-\left[  d(n-1)+\left(  e+1\right)
/2\right]  $ and $\mathcal{D}$ is any constant coefficient differential
operator, then we have $\mathcal{D}f\ \in\mathcal{L}^{2}(\frak{n},d\lambda)$. \label{=corr-growth}
\end{corollary}

\begin{proof}
The group $G$ acts on $I(t,\varepsilon)$ by right translations. The
corresponding action of the Lie algebra $\frak{g}$ in the noncompact picture
$E(t,\varepsilon)$ is by polynomial coefficient vector fields on $\frak{n}$.
The action of $x\in\frak{n}$ is independent of $\left(  t,\varepsilon\right)
$ and is given simply by the directional derivative in the direction $x$. In
particular, the space $E(t,\varepsilon)$ is \emph{invariant} for the action of
constant coefficient differential operators.

Thus$\ f^{\prime}\equiv\mathcal{D}f$ also belongs to $E(t,\varepsilon)$. Thus
$f^{\prime}$ is the restriction to $N$ of a $\overline{P}$-equivariant smooth
function on $G$. Since $G=\overline{P}K$, any such function is
\emph{determined }by its restriction to $K$. The constant function $1$ on $K$
corresponds to the spherical vector $\Phi_{t}$ in $I(t)$. Thus if $c$ is the
maximum of $\left|  f^{\prime}\right|  $ on $K$, then we have $\left|
f^{\prime}\right|  \leq c\Phi_{t}$, and the corollary follows from the
previous lemma.
\end{proof}

We can now prove Propositions \ref{=lemma-phi0} and \ref{=lemma-dmphi} (for
$G\neq Sp_{n}\left(  \mathbb{C}\right)  $.

\begin{proof}
(of Propositions \ref{=lemma-phi0} and \ref{=lemma-dmphi} for $G\neq
Sp_{n}\left(  \mathbb{C}\right)  $) From the table in Section \ref{sec-table}
we see that in every case except $G=Sp_{n}(\mathbb{C})$, we have $2d>e+1$.
Consequently we get
\[
-dn<-\left[  d\left(  n-1\right)  +\left(  e+1\right)  /2\right]  \text{.}%
\]
Proposition \ref{=lemma-phi0} now follows from Lemma \ref{lemma-intcond}, and
Proposition \ref{=lemma-dmphi} follows immediately from Corollary
\ref{=corr-growth} for all groups \emph{except }for $G=Sp_{n}(\mathbb{C})$.
\end{proof}

Suppose now that $G$ is $Sp_{n}(\mathbb{C)}$. Then $L=GL_{n}(\mathbb{C)}$ and
$\frak{n}$ is the space of $n\times n$ complex symmetric matrices. We write
$\mathcal{V}$ for the finite-dimensional space of holomorphic polynomials on
$\frak{n}$ spanned by all the minors of the symmetric matrix $x$, and let
$\varepsilon$ be the unitary character of $L$ given by $\varepsilon
(l)=\dfrac{\det l}{\left|  \det l\right|  }$.

\begin{lemma}
$\mathcal{V}$ is a $Sp_{n}(\mathbb{C})$-invariant subspace of $E(1,\varepsilon
)$.
\end{lemma}

\begin{proof}
The character $e^{\nu}$ of $L$ is simply $\left|  \det l\right|  $. Therefore,
the space $I(1,\varepsilon)$ consists of smooth functions on $G=Sp_{n}%
(\mathbb{C})$ satisfying
\[
f(l\overline{n}g)=\det(l)f(g)\text{.}%
\]

An easy calculation shows that the representation of $G$ on $I(1,\varepsilon)
$ can be expressed in the noncompact picture $E(1,\varepsilon)$ as follows:
\[
g\cdot f(x)=\det(a+xc)\,f\left(  \left[  a+xc\right]  ^{-1}\left[
b+xd\right]  \right)  \text{ for }g=\left[
\begin{array}
[c]{cc}%
a & b\\
c & d
\end{array}
\right]  \in G.
\]

For $p=\left[
\begin{array}
[c]{cc}%
a^{-1} & b\\
0 & a^{t}%
\end{array}
\right]  \in P$ and $w=\left[
\begin{array}
[c]{cc}%
0 & -1\\
1 & 0
\end{array}
\right]  $ we have
\begin{align*}
p\cdot f(x)  &  =\det(a^{-1})\,f\left(  ab+axa^{t}\right) \\
w\cdot f(x)  &  =\det(x)\,f\left(  -x^{-1}\right)
\end{align*}

Evidently, transformations of the form $x\longmapsto c+axa^{t}$ take minors of
$x$ to linear combinations of (possibly smaller) minors; thus $\mathcal{V} $
is $P$-invariant.

Also each minor of $x^{-1}$ is equal to $\pm\det(x)^{-1}$ times the
complementary minor of $x$; thus $\mathcal{V}$ is $w$-invariant. Since $P$ is
a maximal parabolic subgroup, $w$ and $P$ generate $G$, and hence the space
$\mathcal{V}$ is $G$-invariant.

It remains only to prove that $\mathcal{V}\subset E(1,\varepsilon)$. Using the
Gelfand-Naimark decomposition, the functions in $\mathcal{V}$ can be lifted to
$\overline{P}$-equivariant functions on the dense open set $\overline{P}N$ in
$G$. From the $G$-invariance of the finite-dimensional space $\mathcal{V}$, it
follows that these functions transform finitely under right translations by
$K$. Therefore they extend to smooth functions on $K$, and hence on $G$. The
lemma follows.
\end{proof}

\begin{corollary}
For $G=Sp_{n}(\mathbb{C})$, the function $\det(x)$ belongs to the space
$E(1,\varepsilon)$. \endproof
\end{corollary}

We can now finish the proof of Lemma \ref{=lemma-dmphi}.

\begin{proof}
(of Lemma \ref{=lemma-dmphi} for $G=Sp_{n}(\mathbb{C})$) For $G=Sp_{n}%
(\mathbb{C})$, \ we have $d=1$, and the function $\Psi$ is given explicitly
by
\[
\Psi(x)=\det(1+x\overline{x})^{-n/2}%
\]
and we have
\begin{align*}
\phi(x)  &  =\det(x)\det(\overline{x})\\
D  &  =\det(\partial_{x})\det(\partial_{\overline{x}})\text{.}%
\end{align*}
Thus
\begin{equation}
D\Psi=\det(\partial_{x})\det(\partial_{\overline{x}})\det(1+x\overline
{x})^{-n/2}\text{.} \label{=d-phi0}%
\end{equation}

To calculate this, we proceed as follows. First of all, it is well known that
for $u$ a complex symmetric matrix
\[
\det(\partial_{u})\det(u)^{s}=\operatorname{const}\cdot\det(u)^{s-1}%
\]
where the constant can be calculated using, for example, the Capelli identity
from \cite{kostant-sahi}.

Making a simple change of variables, we deduce
\[
\det(\partial_{u})\det(1+u)^{s}=\operatorname{const}\cdot\det(1+u)^{s-1}%
\text{.}%
\]
Now if $v$ is a fixed $n\times n$ complex matrix, then changing variables from
$u$ to $v^{t}uv$, we get
\[
\det(\partial_{u})\det(1+v^{t}uv)^{s}=\operatorname{const}\cdot\det
(vv^{t})\det(1+v^{t}uv)^{s-1}\text{.}%
\]
This can be rewritten as
\[
\det(\partial_{u})\det(1+vv^{t}u)^{s}=\operatorname{const}\cdot\det
(vv^{t})\det(1+vv^{t}u)^{s-1}\text{.}%
\]
By analytic continuation, we get for all complex symmetric $w$
\[
\det(\partial_{u})\det(1+wu)^{s}=\operatorname{const}\cdot\det(w)\det
(1+wu)^{s-1}\text{.}%
\]

Applying this to (\ref{=d-phi0}), we obtain
\[
D\Psi=\operatorname{const}\cdot\det(\partial_{x})\det(x)\det(1+x\overline
{x})^{-n/2-1}.
\]

The function $\det(1+x\overline{x})^{-n/2-1}$ is the spherical vector in
$E(-n-2)$. Also, by the corollary above, $\det(x)$ belongs to $E(1,\varepsilon
)$. Each of these functions extends to a smooth function on $G$ with
appropriate $\overline{P}$-equivariance. By considering the equivariance of
the product, we deduce
\[
\det(x)\det(1+x\overline{x})^{-n/2-1}\in E(-n-1,\varepsilon).
\]

Since $D^{m-1}\det(\partial_{x})$ is a constant coefficient differential
operator, we get
\[
D^{m}\Psi=\left[  \operatorname{const}\cdot D^{m-1}\det(\partial_{x})\right]
\left[  \det(x)\det(1+x\overline{x})^{-n/2-1}\right]  \in E(-n-1,\varepsilon
).
\]
Now in the present case we have $d=1,e=1$, thus we get
\[
-\left[  d\left(  n-1\right)  +\left(  e+1\right)  /2\right]  =-n>-n-1
\]
and so the result follows from Corollary \ref{=corr-growth}.
\end{proof}

\subsection{Proof of the $\mathcal{L}^{1}$ estimate}

We fix $k$ and denote the spherical vector $\Phi_{k,n}=\Phi_{-dk}$ by simply
$\Phi$ as before. In order to prove the necessary estimates for $\Phi$, we
first relate it to the ``rank $1$'' spherical vector
\[
\Upsilon=\Phi_{1,n}=\Phi_{-d}.
\]
We now describe the key result in \cite[Theorem 0.1]{hilbert-two} concerning
the function $\Upsilon$. Let $\tau=\frac{d-e-1}{2}$ and let $K_{\tau}$ be the
corresponding one-variable $K$-Bessel function; define an $M$-invariant
function $\upsilon$ on the rank $1$ orbit $\mathcal{O}_{1}=L\cdot y_{1}%
\subset\overline{\frak{n}}\approx\frak{n}^{\ast}$ by the formula
\[
\upsilon\left(  z\left[  m\cdot y_{1}\right]  \right)  =\frac{K_{\tau}%
(z)}{z^{\tau}}\text{ for }z\in\mathbb{R}^{+}\text{, }m\in M\text{.}%
\]
Then writing $d\mu_{1}$ for the equivariant measure on $\mathcal{O}_{1}$, we
have
\[
\widehat{\upsilon d\mu_{1}}=\Upsilon d\lambda
\]
where $d\lambda$ is the Lebesgue measure on $\frak{n}$, and $\widehat{}$
\ denotes the Fourier transform of tempered distributions. This result is
proved in Propositions 2.1 and 2.2 of \cite{hilbert-two}. For our present
purposes, it is crucial that $\tau$ depends only on $d$ and $e$ but does
\emph{not} depend on $n$.

An immediate consequence of \ref{=lemma-growth} is the relation
\begin{equation}
\Phi=\Upsilon^{k}. \label{=uk}%
\end{equation}
This in turn implies a relation between the Fourier transforms of $\Phi$ and
$\Upsilon$ which we now explain. We start with the following abstract situation:

Suppose $A$ is a Lie group, $\chi$ is a positive character of $A$, and
$B\supset C$ are subgroups such that each of the homogeneous spaces $A/B$ and
$A/C$ admit $\chi$-equivariant measures $dm_{A/B}$ and $dm_{A/C}.$

\begin{lemma}
\label{lemma-Cint} In the above situation, the space $Z=B/C$ admits a
$B$-invariant measure $dz$. Moreover, the formula
\[
\mathcal{C}f\left(  aB\right)  =\int_{Z}f\left(  az\right)  dz
\]
gives a well-defined operator $\mathcal{C}=\mathcal{C}_{A,B,C}:\mathcal{L}%
^{1}\left(  A/C\right)  \rightarrow\mathcal{L}^{1}\left(  A/B\right)  $
satisfying
\begin{equation}
\int_{A/B}\left[  \mathcal{C}f\right]  dm_{A/B}=\int_{A/C}fdm_{A/C}.
\label{=Cint-int}%
\end{equation}
\end{lemma}

\begin{proof}
By Lemma \ref{lemma-equiv} for $b\in B$, $c\in C$, we get
\[
\chi\left(  b\right)  ^{-1}=\frac{\left|  \det\nolimits_{\frak{a}}\left(
\operatorname*{Ad}b\right)  \right|  }{\left|  \det\nolimits_{\frak{b}}\left(
\operatorname*{Ad}b\right)  \right|  }\text{, }\chi\left(  c\right)
^{-1}=\frac{\left|  \det\nolimits_{\frak{a}}\left(  \operatorname*{Ad}%
c\right)  \right|  }{\left|  \det\nolimits_{\frak{c}}\left(
\operatorname*{Ad}c\right)  \right|  }.
\]
Specializing to $b=c$, this implies
\[
\frac{\left|  \det\nolimits_{\frak{b}}\left(  \operatorname*{Ad}c\right)
\right|  }{\left|  \det\nolimits_{\frak{c}}\left(  \operatorname*{Ad}c\right)
\right|  }=1
\]
and another application of \ref{lemma-equiv} proves the existence of an
invariant measure on $B/C$.

The left side of formula (\ref{=Cint-int}) gives a $\chi$-equivariant mean on
the space $C_{c}\left(  A/C\right)  $ thus it agrees with the right side after
suitable normalization of the various measures involved. \ On the other hand
we have
\[
\left|  \mathcal{C}f\left(  a\right)  \right|  =\left|  \int_{Z}f\left(
az\right)  dz\right|  \leq\int_{Z}\left|  f\left(  az\right)  \right|
dz=\mathcal{C}\left|  f\right|  \left(  a\right)
\]

therefore for $f$ in $C_{c}\left(  A/C\right)  $ we get
\[
\int_{A/B}\left|  \mathcal{C}f\right|  dm_{A/B}\leq\int_{A/B}\mathcal{C}%
\left|  f\right|  dm_{A/B}=\int_{A/C}\left|  f\right|  dm_{A/C}.
\]

This shows that $\mathcal{C}$ extends to a bounded linear operator from
$\mathcal{L}^{1}\left(  A/C\right)  \rightarrow\mathcal{L}^{1}\left(
A/B\right)  $ such that the formula (\ref{=Cint-int}) continues to hold.
\end{proof}

We apply the previous result to the situation where
\[
A=L,B=S=\text{stab }\mathbf{y}^{1}\text{, }C=S^{\prime}=\text{stab }%
\mathbf{y}^{\prime}\text{ }%
\]
with $\mathbf{y}^{1}=y_{1}+y_{2}+\ldots+y_{k}$ as before, and $\mathbf{y}%
^{\prime}=(y_{1},y_{2},\ldots,y_{k})\in$ $\mathcal{O}_{1}\times\ldots
\times\mathcal{O}_{1}$. The space $\mathcal{O}=L/S$ is the rank $k$ orbit and
hence by Lemma \ref{lemma-orbitequiv} carries a $e^{2dk\nu}$-equivariant
measure. On the other hand, the space $\mathcal{O}_{1}\times\ldots
\times\mathcal{O}_{1}$ also carries a $e^{2dk\nu}$-equivariant measure, viz.
$d\mu^{\prime}=d\mu_{1}\times\ldots\times d\mu_{1}$; moreover in this
situation $\mathcal{O}^{\prime}=L/S^{\prime}$ is an open subset whose
complement has measure $0$. Thus $\mathcal{O}^{\prime}$ also admits an
$e^{2dk\nu}$-equivariant measure. Thus by the previous lemma, obtain a well
defined operator $\mathcal{C}=\mathcal{C}_{L,S,S^{\prime}}:\mathcal{L}%
^{1}\left(  \mathcal{O}^{\prime}\right)  \rightarrow\mathcal{L}^{1}\left(
\mathcal{O}\right)  $ satisfying formula (\ref{=Cint-int}).

Now given a function $f$ on $\mathcal{O}_{1}$, we define a function $\breve
{f}$ on $\mathcal{O}$ by the following two-step procedure: first define
$\overline{f}$ on $\mathcal{O}^{\prime}$ by
\[
\overline{f}\left(  l\cdot y^{\prime}\right)  =f\left(  l\cdot y_{1}\right)
\cdots f\left(  l\cdot y_{k}\right)  ,
\]
and then set
\[
\breve{f}=\mathcal{C}\overline{f}.
\]

Then we have the following result:

\begin{lemma}
\label{lemma-gdef} For $\upsilon$ as above, put $g=\breve{\upsilon
}=\mathcal{C}\overline{\upsilon}$, then we have
\[
\widehat{gd\mu}=\Phi d\lambda.
\]
\end{lemma}

\begin{proof}
It suffices to prove
\[
\int_{y\in\mathcal{O}}e^{-i\left\langle x,y\right\rangle }g\left(  y\right)
d\mu\left(  y\right)  =\Phi\left(  x\right)  .
\]
To show this, we rewrite the left side as
\begin{equation}
\int e^{-i\left\langle x,l\cdot y^{1}\right\rangle }\mathcal{C}\overline
{\upsilon}\left(  l\cdot y^{1}\right)  d\mu\left(  l\cdot y^{1}\right)  .
\label{left-rewrite}%
\end{equation}
Now we have
\[
\left\langle x,l\cdot y^{1}\right\rangle =\left\langle x,l\cdot y_{1}%
\right\rangle +\cdots+\left\langle x,l\cdot y_{k}\right\rangle .
\]
Thus setting
\[
\eta\left(  l\cdot y_{1}\right)  =\exp\left(  -i\left\langle x,l\cdot
y_{1}\right\rangle \right)  ,
\]
formula (\ref{left-rewrite}) becomes
\[
\int\mathcal{C}\overline{\eta\upsilon}d\mu.
\]
Now using the previous lemma, we can rewrite this as
\[
\int\overline{\eta\upsilon}d\mu^{\prime}=\prod_{j=1}^{k}\left[  \int f\left(
l\cdot y_{j}\right)  \eta\left(  l\cdot y_{j}\right)  d\mu_{1}\left(  l\cdot
y_{j}\right)  \right]  =\Upsilon^{k}=\Phi.
\]
\end{proof}

\begin{proof}
(of Proposition \ref{prop-L1}) In view of the previous lemma, it remains only
to prove that $g\in\mathcal{L}^{1}(\mathcal{O},d\mu)$. In turn, using Lemma
\ref{lemma-Cint}, it suffices to show that $\overline{\upsilon}\in
\mathcal{L}^{1}(\mathcal{O}^{\prime},d\mu^{\prime})$, or equivalently that
\[
\upsilon\in\mathcal{L}^{1}(\mathcal{O}_{1},d\mu_{1}).
\]

This is essentially contained in Proposition 2.1 of \cite{hilbert-two}. The
key point is that by Proposition \ref{prop-sing} for $k=1$, we get
\[
\int_{\mathcal{O}_{1}}\upsilon d\mu_{1}=\int_{\mathbb{R}_{+}}\frac{K_{\tau
}(z)}{z^{\tau}}z^{dn-1}dz.
\]
Since the function $K_{\tau}(z)$ has exponential decay at infinity, it
suffices to prove that the integral on the right converges at $0$. For this we
note that $\frac{K_{\tau}(z)}{z^{\tau}}$ has a pole of order $2\tau=d-e-1$ at
$0$ if $\tau>0$, and a logarithmic singularity if $\tau=0$. At any rate
$\left(  dn-1\right)  -2\tau=d\left(  n-1\right)  +e\ $is greater than $-1$,
which guarantees the convergence of the integral.
\end{proof}

\subsection{Proof of the $\mathcal{L}^{2}$ estimate}

The key to the proof of Proposition \ref{=prop-L2} is a ``stability'' result
for the function $g$ defined in Lemma \ref{lemma-gdef}. To state this, we
temporarily write $g_{k,n}$ and $d\mu_{k,n}$ for $g$ and $d\mu$, in order to
emphasize dependence on $k$ (the rank of the orbit) and $n$ (the rank of the
Jordan algebra). Thus Lemma \ref{lemma-gdef} becomes
\[
\widehat{g_{k,n}d\mu_{k,n}}=\Phi_{k,n}d\lambda.
\]

We now recall the notation $\overline{\frak{n}}_{1},\overline{\frak{n}}_{0}$,
$L_{1}$, $M_{1}$ etc., introduced in subsection \ref{sec-equiv}. Thus
$\overline{\frak{n}}_{1}$ is a Jordan algebra of rank $k$ (with same values of
$d$ and $e$ as $\overline{\frak{n}}$). By applying the considerations of the
previous sections to $\overline{\frak{n}}_{1}$we obtain a family of functions
$g_{j,k};$ $j=1,...,k,$ defined on the various $L_{1}$-orbits in
$\overline{\frak{n}}_{1}$. We are particularly interested in the function
\[
\widetilde{g}=g_{k,k}%
\]
which is defined on the \emph{open }orbit $\widetilde{\mathcal{O}}$ in
$\overline{\frak{n}}_{1}$. Now by definition we have $\overline{\frak{n}}%
_{1}\subset\overline{\frak{n}},$ and moreover we have $\widetilde{\mathcal{O}%
}\subset\mathcal{O}$, where $\mathcal{O}$ is the rank $k$ orbit in
$\overline{\frak{n}}$. Thus we can restrict the function $g=g_{k,n}$ from
$\mathcal{O}$ to $\widetilde{\mathcal{O}}$. The crucial ``stability'' result
is the following:

\begin{lemma}
\label{lemma-stab}With the above notation, we have $g|_{\widetilde
{\mathcal{O}}}=\widetilde{g}$.
\end{lemma}

\begin{proof}
The function $\widetilde{g}$ is also defined by the analogous two-step
procedure applied to the Jordan algebra $\overline{\frak{n}}_{1}$. We start
with the $M_{1}$-invariant function $\widetilde{\upsilon}$ on the rank $1$
orbit $\widetilde{\mathcal{O}}_{1}\subset\overline{\frak{n}}_{1}$
corresponding to the Bessel function $K_{\tau}/z^{\tau}$. As observed after
the definition $\upsilon$, the parameter $\tau=\left(  d-e-1\right)  /2$ is
independent of $n$. Thus we get
\begin{equation}
\upsilon|_{\widetilde{\mathcal{O}}_{1}}=\widetilde{\upsilon}, \label{=u-util}%
\end{equation}
which is the rank $1$ version of the present lemma.

Next, we consider the open $L_{1}$-orbit $\widetilde{\mathcal{O}}^{\prime}$ in
$\widetilde{\mathcal{O}}_{1}\times\cdots\times\widetilde{\mathcal{O}}_{1}$,
and define the analogous function $\overline{\widetilde{\upsilon}}$ by the
formula
\[
\overline{\widetilde{\upsilon}}\left(  l\cdot y^{\prime}\right)
=\widetilde{\upsilon}\left(  l\cdot y_{1}\right)  \cdots\widetilde{\upsilon
}\left(  l\cdot y_{k}\right)
\]
for $y^{\prime}=\left(  y_{1},\cdots,y_{k}\right)  \in\widetilde{\mathcal{O}%
}^{\prime}$ and $l$ in $L_{1}$. Comparing this with the definition of
$\overline{\upsilon}$, and using formula (\ref{=u-util}) we deduce
\[
\overline{\upsilon}|_{\widetilde{\mathcal{O}}^{\prime}}=\overline
{\widetilde{\upsilon}}%
\]
Now the functions $g$ and $\widetilde{g}$ are defined by the integrals
\begin{align}
g\left(  l\cdot y^{1}\right)   &  =\int_{Z}\overline{\upsilon}\left(  l\cdot
z\right)  dz\text{ for }l\text{ in }L\label{g-gtil}\\
\widetilde{g}\left(  l\cdot y^{1}\right)   &  =\int_{\widetilde{Z}}%
\overline{\widetilde{\upsilon}}\left(  l\cdot\widetilde{z}\right)
d\widetilde{z}\text{ for }l\text{ in }L_{1}\nonumber
\end{align}
where $dz$ and $d\widetilde{z}$ are the invariant measures on the homogeneous
spaces $Z=S/S^{\prime}\subset L/S^{\prime}=\mathcal{O}^{\prime}$ and
$\widetilde{Z}=\left(  S\cap L_{1}\right)  /\left(  S^{\prime}\cap
L_{1}\right)  \subset L_{1}/\left(  S^{\prime}\cap L_{1}\right)
=\widetilde{\mathcal{O}}^{\prime}$. However, as in formula (\ref{=stabilizer})
we see that
\begin{align*}
S  &  =\left(  \left(  S\cap L_{1}\right)  \times L_{0}\right)  \cdot
\overline{U,}\\
S^{\prime}  &  =\left(  \left(  S^{\prime}\cap L_{1}\right)  \times
L_{0}\right)  \cdot\overline{U}.
\end{align*}
Thus in the imbedding $\widetilde{\mathcal{O}}^{\prime}\subset\mathcal{O}%
^{\prime}$, we have
\[
Z=\widetilde{Z}.
\]
Moreover, since both measures are $L_{1}$-invariant, we have
\[
dz=d\widetilde{z}.
\]
Thus the two integrals in formula (\ref{g-gtil}) \emph{coincide }for $l$ in
$L_{1}$, and the result follows.
\end{proof}

Let $f\mapsto\check{f}$ denote the inverse Fourier transform which maps
functions on $\frak{n}_{1}$ to functions on $\overline{\frak{n}}_{1}$. Thus
\[
\check{f}\left(  y\right)  =\int_{\frak{n}_{1}}e^{i\left\langle
x,y\right\rangle }f\left(  x\right)  d\lambda
\]
where $d\lambda$ is the Lebesgue measure on $\frak{n}_{1}$.

\begin{lemma}
Writing $\widetilde{\phi}$ for the Jordan norm polynomial on $\overline
{\frak{n}}_{1}$, we have
\[
\left|  (\Phi_{k,k}\check{)}\right|  =\widetilde{g}\left|  \widetilde{\phi
}\right|  ^{d-\left(  e+1\right)  }=\left|  \widetilde{g}\widetilde{\phi
}^{d-\left(  e+1\right)  }\right|  .
\]
\end{lemma}

\begin{proof}
The Fourier transform of tempered distributions is defined by adjointness from
its action on Schwartz functions, and we have the relation
\[
\widehat{\check{f}d\lambda}=fd\lambda.
\]

Now by the definition of $\widetilde{g}$ we have
\[
\widehat{\widetilde{g}d\mu}=\Phi_{k,k}d\lambda,
\]
where $d\mu$ is the equivariant measure on the open orbit $\widetilde
{\mathcal{O}}\subset\overline{\frak{n}}_{1}$. Propositions \ref{prop-open} and
\ref{prop-sing} imply that in polar coordinates, the measures $d\lambda$ and
$d\mu$ are given by $P_{k}^{e+1}V_{k}^{d}d_{k}^{\times}z$ and $P_{k}^{d}%
V_{k}^{d}d_{k}^{\times}z$ respectively. Thus, writing $\widetilde{\phi}$ for
the Jordan norm polynomial on $\overline{\frak{n}}_{1}$, we get
\[
d\mu=\left|  \widetilde{\phi}\right|  ^{d-\left(  e+1\right)  }d\lambda.
\]
Combining these formulas we obtain the result.
\end{proof}

\begin{lemma}
Let $\widetilde{D}$ be the Cayley operator on $\frak{n}_{1}$ then for $l\geq0$
we have
\[
\left|  (\widetilde{D}^{l}\Phi_{k,k}\check{)}\right|  =\left|  \widetilde
{g}\widetilde{\phi}^{l+d-\left(  e+1\right)  }\right|  .
\]
\end{lemma}

\begin{proof}
If $f(x)$ is a Schwartz function on $\frak{n}_{1}$ and $q(y)$ is a homogeneous
polynomial on $\overline{\frak{n}}_{1}$, then we have (up to a scalar
multiple)
\[
(\partial_{q}f\check{)}=q\check{f}\text{,}%
\]
where $\partial_{q}$ is the constant coefficient differential operator on
$\frak{n}_{1}$ with ``symbol'' $q$. \ Thus the proof of the Lemma consists in
establishing that the above identity continues to hold when $f$ \ (like
$\Phi_{k,k}$) is a smooth function of polynomial growth such that $\check{f}$
$\in\mathcal{L}^{1}\left(  \overline{\frak{n}}_{1},d\lambda\right)  $. This is
fairly standard; indeed by adjointness we have the result
\[
\widehat{q\check{f}d\lambda}=\partial_{q}(fd\lambda),
\]
where the derivative on the right is the distributional derivative. Under the
assumption on $f$ , the right side equals $\left(  \partial_{q}f\right)
\,d\lambda$, and the result follows.
\end{proof}

We are now in a position to prove Proposition \ref{=prop-L2}.

\begin{proof}
(of Proposition \ref{=prop-L2}) The function $g$ is $M$-invariant, thus by
Proposition \ref{prop-sing}, it suffices to prove the convergence of the
integral
\[
\int_{\mathcal{C}_{k}}\left|  g(z_{1}y_{1}+\ldots+z_{k}y_{k})\right|
^{2}\left[  P_{k}^{n-k+1}V_{k}\right]  ^{d}d_{k}^{\times}z.
\]
By the previous Lemma, this can be rewritten as
\[
\int_{\mathcal{C}_{k}}\left|  \widetilde{g}(z_{1}y_{1}+\ldots+z_{k}%
y_{k})\right|  ^{2}\left[  P_{k}^{n-k+1}V_{k}\right]  ^{d}d_{k}^{\times}z.
\]
Using Proposition \ref{prop-sing} we can further rewrite this as
\[
\int_{\overline{\frak{n}}_{1}}\left|  \widetilde{g}\right|  ^{2}\left|
\widetilde{\phi}^{t}\right|  d\lambda\text{ where }t=d(n-k+1)-\left(
e+1\right)  \text{.}%
\]

Thus it suffices to prove that
\begin{equation}
\left|  \widetilde{g}^{2}\widetilde{\phi}^{t}\right|  \in\mathcal{L}%
^{1}\left(  \overline{\frak{n}}_{1},d\lambda\right)  . \label{t-L1}%
\end{equation}

Now the map $f\mapsto\check{f}$ extends as an isometry from $\mathcal{L}%
^{2}\left(  \frak{n}_{1},d\lambda\right)  $ to $\mathcal{L}^{2}\left(
\overline{\frak{n}}_{1},d\lambda\right)  $ (after suitable normalizations of
the Lebesgue measures). Thus we have
\begin{equation}
f_{1},f_{2}\in\mathcal{L}^{2}\left(  \frak{n}_{1},d\lambda\right)
\Rightarrow\check{f}_{1}\check{f}_{2}\in\mathcal{L}^{1}\left(  \overline
{\frak{n}}_{1},d\lambda\right)  \text{;} \label{f1f2L1}%
\end{equation}
we shall deduce (\ref{t-L1}) from (\ref{f1f2L1}) by a suitable choice of
$f_{1},f_{2}$.

Let us put
\[
s=t+2\left(  e+1-d\right)  =d(n-k-1)+\left(  e+1\right)  \text{;}%
\]
since $n>k$, we have $s>0$. Now if we set
\begin{equation}
f_{1}=\widetilde{D}^{l_{1}}\Phi_{k,k}\text{, }f_{2}=\widetilde{D}^{l_{2}}%
\Phi_{k,k}\text{ where }l_{1}+l_{2}=s; \label{f1f2}%
\end{equation}
then by the previous lemma we have
\begin{equation}
\left|  \check{f}_{1}\check{f}_{2}\right|  =\widetilde{g}^{2}\left|
\widetilde{\phi}\right|  ^{l_{1}+l_{2}-2(e+1-d)}=\left|  \widetilde{g}%
^{2}\widetilde{\phi}^{t}\right|  \text{.} \label{f1f2t}%
\end{equation}

We now consider two cases: if $G\neq Sp_{n}(\mathbb{C})$ we set $l_{1}=0$ and
$l_{2}=s$; if $G=$ $Sp_{n}(\mathbb{C})$, we set $l_{1}=1$ and $l_{2}=s-1$. In
the former case we have $l_{1},l_{2}\geq0$; while in the latter case, we have
$e=1$, whence $s\geq2$ and $l_{1},l_{2}\geq1$. Thus in either case by the open
orbit estimates of Lemmas \ref{=lemma-phi0} and \ref{=lemma-dmphi}, applied to
the Jordan algebra $\frak{n}_{1}$, we deduce that the functions $f_{1}$ and
$f_{2}$ from formula (\ref{f1f2}) belong to $\mathcal{L}^{2}(\frak{n}%
_{1},d_{1}x)$. Thus formula (\ref{t-L1}) follows from (\ref{f1f2t}) and
(\ref{f1f2L1}).
\end{proof}

\section{Proof of the main results}

We now explain how to deduce Theorems \ref{theoremA} and \ref{theoremB} from
the previous results. As explained in the introduction, the arguments are very
similar to those in \cite{sahi-expl}, \cite{tens} and \cite{hilbert-two}. Thus
we shall limit ourselves to only sketching the proofs of the various results below.

\subsection{Proof of Theorem \ref{theoremA}}

In order to prove Theorem \ref{theoremA}, we introduce a number of spaces.

First of all, let $I\left(  -dk\right)  \subset C^{\infty}(\frak{n})$ be the
space of smooth vectors in the degenerate principal series defined in section
\ref{sec-dp}.The representation $\pi=\pi_{-dk\nu}$ of the group $G$ on this
space is by ``fractional linear transformations'', and we have
\begin{align*}
\left[  \pi\left(  l\right)  f\right]  \left(  x\right)   &  =e^{-dk\nu
}\left(  l\right)  f\left(  \operatorname*{Ad}l^{-1}\left[  x\right]  \right)
\text{ for }l\text{ in }L,\\
\left[  \pi\left(  \exp x^{\prime}\right)  f\right]  \left(  x\right)   &
=f\left(  x+x^{\prime}\right)  \text{ for }x^{\prime}\text{ in }\frak{n.}%
\end{align*}
By \cite{sahi-dp}, the space $E(-dk)$ has an irreducible \emph{unitarizable}
spherical $(\frak{g},K)$-submodule $V$ which we also regard as a subspace of
$C^{\infty}(\frak{n})$. Thus by Harish-Chandra theory, the Hilbert space
closure $\mathcal{H}$ of $V$ with respect to the $(\frak{g},K)$-invariant norm
carries an irreducible unitary representation of $G$. \ 

For convenience, we first describe $\mathcal{H}\ $as the closure of a
$G$-invariant space. For this we introduce the space $\mathbf{V}$ consisting
of those vectors in $I(-dk)$ whose restriction to $K$, and subsequent
expansion in $K$-isotypic components only involves the $K$-types of $V.$ Since
$V$ is $(\frak{g},K)$-invariant, the space $\mathbf{V}$ is $G$-invariant and
we have the following result.

\begin{lemma}
The functions in $\mathbf{V}$ have finite $\mathcal{H}$-norm and
$\mathcal{H}\ $is the closure of $\mathbf{V}$.
\end{lemma}

\begin{proof}
This is a consequence of a general result due to Casselman-Wallach on the
smooth vectors of a representation. In the present situation, one can also
give an alternative proof along the lines of the remark in section 2.4 of
\cite{hilbert-two} as follows.

First of all, the $K$-types of $V$ have multiplicity $1$, and have highest
weights of the form
\[
m_{1}\gamma_{1}+\cdots+m_{k}\gamma_{k}\text{ ,}%
\]
where $m_{1}\geq\cdots\geq m_{k}\geq0$ and $\gamma_{1},\cdots,\gamma_{k}$ are
as in subsection \ref{sec-KM}. Moreover the $\mathcal{H}$-norm on each
$K$-type is computed explicitly in \cite{sahi-dp} and the ratio of the
$\mathcal{H}$-norm to the $\mathcal{L}^{2}(K)$-norm grows at most polynomially
in $\left(  m_{1},\cdots,m_{k}\right)  $. On the other hand by the
Riemann-Lebesgue lemma for $f$ in $\mathbf{V}$\textbf{,} the $\mathcal{L}%
^{2}(K)$-norms of its $K$-isotypic components decay rapidly. Thus such an $f$
will have finite $\mathcal{H}$-norm. Evidently since $V\subset\mathbf{V}$, the
closure of $\mathbf{V}$ is $\mathcal{H}$ as well.
\end{proof}

Next recall the space $\mathcal{H}_{\mathcal{O}}=\mathcal{L}^{2}%
(\mathcal{O},d\mu)$; by Mackey theory, this space carries a natural
irreducible unitary representation\ $\pi_{\mathcal{O}}$ of $P$, which is given
by the following explicit formulas:
\begin{align*}
\left[  \pi_{\mathcal{O}}\left(  l\right)  \psi\right]  \left(  y\right)   &
=e^{dk\nu}\left(  l\right)  \psi\left(  \operatorname*{Ad}l^{-1}\left[
y\right]  \right)  \text{ for }l\text{ in }L\\
\left[  \pi_{\mathcal{O}}\left(  \exp x\right)  \psi\right]  \left(
y\right)   &  =e^{-i\left\langle x,y\right\rangle }\psi\left(  y\right)
\text{ for }x\text{ in }\frak{n}\text{ ,}%
\end{align*}
where $\left\langle x,y\right\rangle $ is the normalized Killing form of
subsection \ref{sec-LH}. We shall prove Theorem \ref{theoremA} by constructing
a unitary $P$-isomorphism $\mathcal{I}$ between $\left(  \pi|_{P}%
,\mathcal{H}\right)  $ and $\left(  \pi_{\mathcal{O}},\mathcal{H}%
_{\mathcal{O}}\right)  $.

We first define $\mathcal{I}$ on a suitable subspace of $\mathcal{H}$. For
this, let $\mathcal{C}\left(  G\right)  $ be the convolution algebra of smooth
$\mathcal{L}^{1}$ functions on $G$. Then by standard arguments, $\pi$ extends
to a representation of $\mathcal{C}\left(  G\right)  $ on $\mathbf{V}$ and we
define
\[
\mathbf{W}=\pi\left(  \mathcal{C}\left(  G\right)  \right)  \Phi
\subset\mathbf{V}\text{ ,}%
\]
where $\Phi=\Phi_{-dk}$ is the spherical vector in $I\left(  -dk\right)  $.
Since $G=PK$ and $\Phi$ is $K$-fixed, we also have
\[
\mathbf{W}=\pi\left(  \mathcal{C}\left(  P\right)  \right)  \Phi;
\]
and we shall prove the following result:

\begin{lemma}
For each $f$ in $\mathbf{W}$ there is a unique $\mathcal{I}$ $\left(
f\right)  \in$ $\mathcal{H}_{\mathcal{O}}$ such that we have the equality
\[
fd\lambda=\widehat{\mathcal{I}\left(  f\right)  d\mu},
\]
of tempered distributions. Furthermore, for all $F\in\mathcal{C}\left(
P\right)  $ we have
\begin{equation}
\mathcal{I}\circ\pi\left(  F\right)  =\pi_{\mathcal{O}}\left(  F\right)
\circ\mathcal{I}\text{ }. \label{=inter}%
\end{equation}
\end{lemma}

\begin{proof}
The key step is, of course, \ref{=prop-L2} which shows that for the function
$\Phi=\Phi_{-dk}$ we have
\[
\Phi d\lambda=\widehat{\psi d\mu}.
\]
where $\psi\in$ $\mathcal{H}_{\mathcal{O}}$; or, equivalently,
\[
\Phi(x)=\int_{\mathcal{O}}e^{-i\left\langle x,y\right\rangle }\psi(y)d\mu(y).
\]
Now for $l$ in $L$, by Lemma \ref{lemma-orbitequiv} we have
\begin{align*}
\int_{\mathcal{O}}e^{-i\left\langle x,y\right\rangle }\left[  \pi
_{\mathcal{O}}\left(  l\right)  \psi\right]  (y)d\mu(y)  &  =e^{dk\nu}\left(
l\right)  \int_{\mathcal{O}}e^{-i\left\langle x,y\right\rangle }\psi\left(
\operatorname*{Ad}l^{-1}\left[  y\right]  \right)  d\mu(y)\\
&  =e^{dk\nu}\left(  l\right)  e^{-2dk\nu}\left(  l\right)  \int_{\mathcal{O}%
}e^{-i\left\langle x,\operatorname*{Ad}l\left[  y\right]  \right\rangle }%
\psi(y)d\mu(y)\\
&  =e^{-dk\nu}\left(  l\right)  \int_{\mathcal{O}}e^{-i\left\langle
\operatorname*{Ad}l^{-1}\left[  x\right]  ,y\right\rangle }\psi(y)d\mu
(y)=\pi\left(  l\right)  f
\end{align*}
Similarly for $x^{\prime}$ in $\frak{n}$, we have
\[
\int_{\mathcal{O}}e^{-i\left\langle x,y\right\rangle }\left[  \pi
_{\mathcal{O}}\left(  \exp x^{\prime}\right)  \psi\right]  \left(  y\right)
d\mu(y)=\int_{\mathcal{O}}e^{-i\left\langle x+x^{\prime},y\right\rangle }%
\psi(y)d\mu(y)=\pi\left(  \exp x^{\prime}\right)  f.
\]

Thus for any $F\in\mathcal{C}\left(  P\right)  $, we have
\[
\left[  \pi\left(  F\right)  \Phi\right]  d\lambda=(\left[  \pi_{\mathcal{O}%
}\left(  F\right)  \psi\right]  d\mu\widehat{)}%
\]
and we can define $\mathcal{I}$ by the formula
\[
\mathcal{I}\left(  \pi\left(  F\right)  \Phi\right)  =\pi_{\mathcal{O}}\left(
F\right)  \psi.
\]
Then $\mathcal{I}$ satisfies the conditions of the lemma. The uniqueness is clear.
\end{proof}

We can now finish the proof of Theorem \ref{theoremA}.

\begin{proof}
(of Theorem \ref{theoremA}) Given the previous lemma, the proof of the result
proceeds along lines similar to \cite{sahi-expl} and \cite[ ]{hilbert-two}. By
the previous lemma, the space $\mathbf{W}_{1}=\mathcal{I}\left(
\mathbf{W}\right)  $ is a $\mathcal{C}\left(  P\right)  $-invariant subspace
of $\mathcal{H}_{\mathcal{O}}$, and moreover we can equip it with a second
$P$-invariant norm, namely that transferred from $\mathcal{H}$. Now as
explained in \cite[3.3]{sahi-expl}, it follows from \cite{poguntke} that
$\mathbf{W}_{1}$contains a further $\mathcal{C}\left(  P\right)  $-invariant
subspace $\mathbf{W}_{2}$ on which the two norms coincide (up to a scalar
multiple which we normalize to be $1$ by rescaling $\mathcal{I}$).

Since $\mathcal{H}_{\mathcal{O}}$ is irreducible, $\mathbf{W}_{2}$ is dense in
$\mathcal{H}_{\mathcal{O}}$ and thus $\mathcal{H}_{\mathcal{O}}$ can be
regarded as the closure of $\mathbf{W}_{2}$ with respect to the $\mathcal{H}%
$-norm. It follows that the two norms agree on $\mathbf{W}_{1}$ as well, and
thus the map
\[
\mathcal{I}^{-1}:\mathbf{W}_{1}\rightarrow\mathbf{W}%
\]
extend to an isometric $P$-invariant imbedding $\mathcal{J}$ of $\mathcal{H}%
_{\mathcal{O}}$ into $\mathcal{H}$. Now the image of $\mathcal{J}$ \ is
closed, and contains a $G$-invariant subspace (namely $\mathbf{W})$; thus
since $\mathcal{H}$ is an irreducible representation, it follows that
$\mathcal{J}$ is surjective as well. Thus $\mathcal{J}$ is a unitary
intertwining operator between $\left(  \pi_{\mathcal{O}},\mathcal{H}%
_{\mathcal{O}}\right)  $ and $\left(  \pi|_{P},\mathcal{H}\right)  $. The
required extension of $\left(  \pi_{\mathcal{O}},\mathcal{H}_{\mathcal{O}%
}\right)  $ to $G$ is now given by simply transferring the representation from
$\left(  \pi,\mathcal{H}\right)  $ via $\mathcal{J}^{-1}$.
\end{proof}

\subsection{Proof of Theorem \ref{theoremB}}

We now study tensor products of our representations $\pi_{\mathcal{O}}$. The
analogous study for conformal groups of \emph{Euclidean }Jordan algebras was
conducted in \cite{tens}. Since the statements and proofs from \cite{tens} can
be transferred to our present (non-Euclidean) setting without substantial
changes, we will only sketch some of the arguments below.

Fix $s\geq2$ and a collection of positive integers $k_{1},\ldots,k_{s}$
satisfying the condition
\[
k=k_{1}+\cdots+k_{s}\leq n\text{.}%
\]
For each $i=1,\ldots,s$, let $\mathcal{O}^{i}$ be the $L$-orbit on
$\overline{\frak{n}}$ of rank $k_{i}$, with $L$-equivariant measure $d\mu^{i}%
$. Let $\pi_{\mathcal{O}^{i}}$ be the unitary representation of $G$ on the
space $\mathcal{L}^{2}(\mathcal{O}^{i},d\mu^{i})$ as described in theorem
\ref{theoremA}. We wish to study the tensor product representation
\[
\Pi=\pi_{\mathcal{O}^{1}}\otimes\ldots\otimes\pi_{\mathcal{O}^{s}}\text{ }%
\]
which can be realized explicitly on the space $\mathcal{L}^{2}(\mathcal{O}%
^{1}\times\ldots\times\mathcal{O}^{s},d\mu^{1}\times\ldots\times d\mu^{s})$.

Let $y_{1},\cdots,y_{n}$ be as in subsection \ref{sec-Striples}, and define
\[
v_{i}=y_{m_{i}+1}+y_{m_{i}+2}+\ldots+y_{m_{i}+k_{i}}\text{, where }m_{i}%
=k_{1}+\ldots+k_{i-1},1\leq i\leq s.
\]
Then $v_{i}$ is an orbit representative for $\mathcal{O}^{i}$; $v=$
$v_{1}+\cdots+v_{s}$ is an orbit representative for the rank $k$ orbit
$\mathcal{O}$; and the $L$-orbit of
\[
v^{\prime}=\left(  v_{1},\cdots,v_{s}\right)
\]
is an open subset of $\mathcal{O}^{1}\times\ldots\times\mathcal{O}^{s}$ with
full measure. We denote by $S^{\prime}$ and $S$ the isotropy subgroups of
$v^{\prime}$ and $v$, respectively. In the notation of subsection
\ref{sec-equiv}, we have $v=\mathbf{y}^{1}$, and thus
\[
S=(H_{1}\times L_{0})\cdot U\text{.}%
\]
It is easy to see that $S^{\prime}$ can then be written as
\[
S^{\prime}=(H_{1}^{\prime}\times L_{0})\cdot U\text{,}%
\]
where $H_{1}^{\prime}$ is a certain reductive subgroup of $H_{1}$. We now
change the notation slightly and write $G^{\prime}$ for $H_{1}$ and
$H^{\prime}$ for $H_{1}^{\prime}$.

\textbf{Example.} Take $G=E_{7(7)}$, $s=2$ and $k_{1}=1$, $k_{2}=2$. Then
$k=n=3$ and $S=G^{\prime}$ (the stabilizer of the identity element of
$\overline{\frak{n}}$ -- the exceptional Jordan algebra of dimension 27). In
this case we have $G^{\prime}=F_{4(4)}$ and $S^{\prime}=H^{\prime
}=\operatorname{Spin}_{4,5}$ (cf. \cite[p. 119]{adams}).

In general, $X=G^{\prime}/H^{\prime}$ is a reductive homogeneous space, and we
write $\operatorname*{Ind}_{H^{\prime}}^{G^{\prime}}1$ for the quasiregular
representation of $G^{\prime}$ on $\mathcal{L}^{2}(X)$. We decompose this
using the Plancherel measure $d\rho$ and the corresponding multiplicity
function $m:\widehat{H}\rightarrow\{0,1,2,\ldots,\infty\}$, i.e.,
\[
\operatorname*{Ind}\nolimits_{H^{\prime}}^{G^{\prime}}1\simeq\int
_{\widehat{G^{\prime}}}^{\oplus}m(\sigma)\sigma\,d\rho(\sigma)\text{.}%
\]
We define a map $\Theta$ from irreducible unitary representations of
$G^{\prime}$ to unitary representations of$\ P$ defined as follows
\[
\Theta(\sigma)=\operatorname*{Ind}\nolimits_{SN}^{P}(E\sigma\otimes\chi_{v}),
\]
where $E\sigma$ denotes the trivial extension of $\sigma$ to $S=\left(
G^{\prime}\times L_{0}\right)  \cdot U$, and $\chi_{v}$ is the unitary
character of $N$ defined by
\[
\chi_{v}\left(  \exp x\right)  =e^{-i\left\langle v,x\right\rangle }.
\]
An easy application of Mackey theory shows that all representations
$\Theta(\sigma)$ are unitary irreducible representations of $P$, and
$\Theta(\sigma)\simeq\Theta(\sigma^{\prime})$ if and only if $\sigma
\simeq\sigma^{\prime}$.

\begin{proposition}
The restriction of $\Pi$ to $P$ decomposes as follows
\begin{equation}
\Pi|_{P}\simeq\int_{\widehat{G^{\prime}}}^{\oplus}m(\sigma)\Theta
(\sigma)\,d\rho(\sigma)\text{,} \label{=Pi-to-P}%
\end{equation}
\end{proposition}

\begin{proof}
This is proved as in \cite[Lemma 2.1]{tens} --- here is a sketch of the
argument. We define an operator $F$ from the space of $\Pi$ to functions on
$P$ by the formula
\[
\left[  Ff\right]  \left(  ln\right)  =\chi_{v}\left(  lnl^{-1}\right)
f(l\cdot v^{\prime}),\text{ }l\in L,\,n\in N.
\]
It is an easy exercise to verify that $F$ gives a unitary isomorphism
\[
\Pi|_{P}\simeq\operatorname*{Ind}\nolimits_{S^{\prime}N}^{P}(1\otimes\chi
_{v}).
\]

Next, using induction by stages we obtain an isomorphism
\[
\operatorname*{Ind}\nolimits_{S^{\prime}N}^{P}(1\otimes\chi_{v})\simeq
\operatorname*{Ind}\nolimits_{SN}^{P}\left(  (\operatorname*{Ind}%
\nolimits_{S^{\prime}}^{S}1)\otimes\chi_{v}\right)  .
\]
A final easy calculation shows that
\[
\operatorname*{Ind}\nolimits_{S^{\prime}}^{S}1\simeq E\left(
\operatorname*{Ind}\nolimits_{H^{\prime}}^{G^{\prime}}1\right)  \simeq
\int_{\widehat{G^{\prime}}}^{\oplus}m(\sigma)\left(  E\sigma\right)
d\rho(\sigma)
\]
Combining the various isomorphisms, we obtain the result.
\end{proof}

Let $\kappa$ be a unitary representation of $G$ on a Hilbert space
$\mathcal{H}$, and $R$ be a subgroup of $G$. We shall write $\mathcal{A}%
(\kappa,R)$ for the von Neumann algebra generated by the operators
$\{\kappa(g)|g\in R\}$. If $G$ is a type I group \cite{mackey}, then for
irreducible $\kappa$ one has $\mathcal{A}(\kappa,G)=\mathcal{B}\left(
\mathcal{H}\right)  $ --- the full algebra of bounded operators on
$\mathcal{H}$. To extend the $P$-decomposition of $\Pi$ from formula
(\ref{=Pi-to-P}) to the $G$- decomposition, we require the following

\begin{proposition}
$\mathcal{A}(\Pi,G)=\mathcal{A}(\Pi,P)$. \label{=prop-vN}
\end{proposition}

Proposition \ref{=prop-vN} was proved in \cite[4.4]{tens} for conformal groups
of Euclidean Jordan algebras. The proof given in \cite{tens} combines the low
rank theory of \cite{li-singular}, \cite{li-lowrank} for classical groups, and
Jordan algebra techniques for the exceptional groups. The arguments extend to
our present setting without any significant modifications. For the readers
convenience, we outline the steps of the argument in Appendix
\ref{app-lowrank}.

\begin{proof}
(of Theorem \ref{theoremB}) Consider the direct integral decomposition of
$\Pi$
\[
\Pi=\int^{\oplus}m(\kappa)\kappa\,d\eta(\kappa)
\]
into irreducible representations of $G$. Then
\[
\mathcal{A}(\Pi,P)\subseteq\int^{\oplus}m(\kappa)\mathcal{A}(\kappa
,P)\,d\eta(\kappa)\subseteq\int^{\oplus}m(\kappa)\mathcal{A}(\kappa
,G)\,d\eta(\kappa)=\mathcal{A}(\Pi,G)\text{.}%
\]
The equality of Proposition \ref{=prop-vN} is possible only when the following
conditions are satisfied (for almost every $\kappa$ with respect to $d\eta$):

\begin{itemize}
\item $\kappa|_{P}$ is irreducible (then $\mathcal{A}(\kappa,P)\,=\mathcal{A}%
(\kappa,G)$);

\item  If $\kappa|_{P}\simeq\kappa^{\prime}|_{P}$, then $\kappa\simeq
\kappa^{\prime}$.
\end{itemize}

In other words, in this case (almost) every irreducible representation
$\Theta(\sigma)$ from the spectrum of $\Pi|_{P}$ extends \emph{uniquely} to a
certain irreducible representation of $G$, which we denote by $\mathbf{\theta
}(\sigma)$; and the $P$-decomposition (\ref{=Pi-to-P}) gives rise to the
$G$-decomposition
\[
\Pi=\int_{\widehat{G^{\prime}}}^{\oplus}m(\sigma)\mathbf{\theta}%
(\sigma)\,d\rho(\sigma)
\]
and the theorem follows.
\end{proof}

\textbf{Example. }Again, take $G=E_{7(7)}$, $s=2$ and $k_{1}=1$, $k_{2}=2$.
Then the map $\sigma\mapsto$ $\mathbf{\theta}(\sigma)$ establishes a
correspondence between the spectrum of $\Pi$ and the spectrum of the rank $1$
reductive symmetric space $F_{4(4)}/Spin(4,5)$. In other words, we obtain a
duality between (some subsets of) the unitary duals of two exceptional groups:
split $F_{4}$ on one side and split $E_{7}$ on the other side. As with Howe's
duality correspondence (the usual $\theta$-correspondence), we expect that
this new duality will have smooth and global analogues.\label{sec-tensor}

\section{Appendix}

\subsection{Measures and Jacobians}

If $F:X\rightarrow Y$ is a diffeomorphism between manifolds and $dx$ is a
measure on $X$, the pushforward of $dx$ is the measure $F_{\ast}dx$ on $Y$
defined by
\[
\int_{Y}f\left[  F_{\ast}dx\right]  =\int_{X}f\left(  F\left(  x\right)
\right)  dx.
\]
If $d\lambda$ is the Lebesgue measure on $\mathbb{R}^{n}$ and $A:U\rightarrow
V$ is a diffeomorphism between open sets in $\mathbb{R}^{n}$, then the
``change of variables'' formula for the Lebesgue integral says
\[
\int_{U}f\left(  A\left(  x\right)  \right)  \left|  \det D_{A}\left(
x\right)  \right|  d\lambda=\int_{V}fd\lambda
\]
where $D_{A}\left(  x\right)  $ $:\mathbb{R}^{n}\rightarrow\mathbb{R}^{n}$ is
the differential of $A$ at $x$ in $U$. In other words, we have
\[
A_{\ast}\left[  \left|  \det D_{A}\left(  x\right)  \right|  d\lambda\right]
=d\lambda\text{ or }%
\]
or equivalently
\[
A_{\ast}d\lambda=gd\lambda\text{ where }g\left(  A\left(  x\right)  \right)
=\left|  \det D_{A}\left(  x\right)  \right|  ^{-1}%
\]

A measure on an open set $U\subset$ $\mathbb{R}^{n}$ will be called
\emph{regular }if it is of the form $\phi d\lambda$\ where $\phi$ is a smooth
positive function. If $\ A:U\rightarrow V$ is a diffeomorphism, then we have
\[
A_{\ast}\left(  \phi d\lambda\right)  =\psi d\lambda\text{ where }\psi\left(
A\left(  x\right)  \right)  =\phi\left(  x\right)  /\left|  \det D_{A}\left(
x\right)  \right|  ;
\]
thus $A_{\ast}$ maps regular measures to regular measures. If $dx$ is a
measure on a smooth manifold $X$\ and $U$ is a coordinate open set, then the
pushforward of $dx|_{U}$ under the coordinate map is a measure on an open set
in $\mathbb{R}^{n}$, which we shall call the ``local expression'' of $dx$. We
shall say that $dx$ is regular if each of these local expressions is a regular
measure in the earlier sense. It follows from our discussion that if
$F:X\rightarrow Y$ is a diffeomorphism and $dx$ is a regular measure on $X$,
then $F_{\ast}dx$ is a regular measure on $Y$.

Suppose $X$ and $Y$ are manifolds with regular measures $dx$ and $dy$, and
$F:$ $X$ $\rightarrow Y$ is a diffeomorphism. The \emph{Jacobian} of $F$ is
the function $J_{F}\left(  x\right)  =J_{F,dx,dy}\left(  x\right)  $ on $X$
satisfying
\[
\int f\left(  y\right)  dy=\int J_{F}\left(  x\right)  f\left(  F\left(
x\right)  \right)  dx.
\]
or, equivalently
\[
F_{\ast}\left[  J_{F}dx\right]  =dy
\]
If we have another diffeomorphism $G:$ $Y$ $\rightarrow Z$, $\ $where $Z$ is a
manifold with regular measure $dz$, then
\begin{equation}
J_{GF}\left(  x\right)  =J_{F}(x)J_{G}\left(  y\right)  \text{ where
}y=F\left(  x\right)  \label{prodjac}%
\end{equation}
If $X$ and $Y$ are open sets in $\mathbb{R}^{n}$, and $dx=$ $\phi
d\lambda,dy=\psi d\lambda$, then it is easy to see that
\[
J_{F,dx,dy}=\frac{\left|  \det D_{F}\left(  x\right)  \right|  \psi\left(
x\right)  }{\phi\left(  F\left(  x\right)  \right)  }.
\]
If $X$ and $Y$ are smooth manifolds, then we can determine the Jacobian in a
similar manner by passing to local coordinates and using formula
(\ref{prodjac}). In particular we see that for regular measures, the Jacobian
is a well-defined smooth positive function.

The following lemma will be quite useful in calculating Jacobians.

\begin{lemma}
\label{lemma-jacdet} Suppose $X$ and $Y$ are manifolds with regular measures
$dx$ and $dy$; let $x\in X$, $y\in Y$, and fix linear bases for the tangent
spaces $T_{x}X$, $T_{y}Y$. Then there is a positive constant $c$, such that
for any diffeomorphism $X\overset{F}{\rightarrow}Y$ satisfying $F\left(
x\right)  =y$, we have
\[
\text{ }J_{F}(x)=c\left|  \det D_{F}\left(  x\right)  \right|
\]
where we regard the differential $D_{F}\left(  x\right)  :$ $T_{x}X\rightarrow
T_{y}Y\ $as a matrix for the above bases.
\end{lemma}

\begin{proof}
If $X$ and $Y$ are open sets in $\mathbb{R}^{n}$, and $dx=$ $\phi
d\lambda,dy=\psi d\lambda$, then
\[
J_{F,dx,dy}=\frac{\left|  \det D_{F}\left(  x\right)  \right|  \psi\left(
x\right)  }{\phi\left(  y\right)  }=c\left|  \det D_{F}\left(  x\right)
\right|  ,
\]
where the determinant is computed for the standard basis of $T_{x}X=T_{y}Y=$
$\mathbb{R}^{n}$. If we use different bases then the scalar $c$ is replaced by
a different scalar, which is still independent of $F$. Passing to local
coordinates, we obtain the result in general.
\end{proof}

Now suppose $X$ is a smooth manifold with the action of a Lie group $G$, and
let $\chi$ be a positive multiplicative character of $G$. A regular measure
$dx$ is called $\chi$-\emph{equivariant}
\[
g_{\ast}dx=\chi\left(  g\right)  dx\text{ for all }g\text{ in }G.
\]
Equivalently
\[
J_{g,dx,dx}\left(  x\right)  =\chi\left(  g\right)  ^{-1}\text{ for all }x\in
X,g\in G.
\]
For example, the Lebesgue measure $d\lambda$ on $\mathbb{R}^{n}$ is
equivariant by the character $\left|  \det\right|  ^{-1}$ of the group
$GL\left(  \mathbb{R}^{n}\right)  $.

If $X$ is a homogeneous space for $G$ then for $x$ in $X$, the tangent space
$T_{x}X$ can be naturally identified with $\frak{g}/\frak{g}^{x}$, where
$\frak{g}^{x}$ is the Lie algebra of the stabilizer $G^{x}$ of $x$. Moreover,
for $g$ in $G$ we have
\[
D_{g}=\operatorname*{Ad}g:\frak{g}/\frak{g}^{x}\rightarrow\frak{g}%
/\frak{g}^{g\cdot x}.
\]

\begin{lemma}
\label{lemma-equiv} A $G$-homogeneous space $X$ admits a $\chi$-equivariant
measure if and only if
\[
\chi\left(  h\right)  ^{-1}=\left|  \det\nolimits_{\frak{g}/\frak{g}^{x}%
}\left(  \operatorname*{Ad}h\right)  \right|  \text{ for all }h\in G^{x}.
\]
\end{lemma}

This result is well-known and can be proved in a manner analogous to Theorem
I.1.9 in \cite{helgason}. Here is a quick argument for the necessity of the
condition. If $h$ is in $G^{x}$, then we have
\[
\chi\left(  h\right)  ^{-1}=J_{h}\left(  x\right)  =c\left|  \det
D_{h}\right|  =c\left|  \det\nolimits_{\frak{g}/\frak{g}^{x}}\left(
\operatorname*{Ad}h\right)  \right|  .
\]
Specializing to $h=1\in G$, we deduce that $c=1$.

\begin{corollary}
\label{cor-equiv} If $G$ is reductive, then the condition of the previous
lemma becomes
\[
\chi\left(  h\right)  =\left|  \det\nolimits_{\frak{g}^{x}}\left(
\operatorname*{Ad}h\right)  \right|  \text{ for all }h\in G^{x}.
\]
\end{corollary}

\subsection{Low rank representations}

Let $\tau$ be a unitary representation of $G.$ Consider its restriction to
$P$, and its further restriction to $N$. Since $N$ is abelian, the restriction
$\tau|_{N}$ decomposes into a direct integral of unitary characters of $N$.
This decomposition defines a projection valued measure on the dual space
$N^{\ast}$, which we identify with $\overline{\frak{n}}$. If this measure is
supported on a single non-open orbit $\mathcal{O}_{r}\subset\overline{N}$, we
say that $\tau$ a \emph{low-rank }representation of $G$ and write
\[
\operatorname*{rank}\nolimits_{N}\tau=r\text{.}%
\]

An element $x_{1}$ is a primitive idempotent in a Jordan algebra $N$, and we
can consider the associated Peirce decomposition
\[
N=N(x_{1},1)+N(x_{1},\frac{1}{2})+N(x_{1},0)\text{.}%
\]
Observe that the spaces $N(x_{1},1)$ and $N(x_{1},0)$ are the Jordan algebras
of ranks $1$ and $n-1$, respectively, with the respect to the Jordan structure
inherited from $N.$

We will write $N_{1}$ and $N_{0}$ for $N(x_{1},1)$ and $N(x_{1},0)$,
respectively. Similarly, we write $G_{0}$ for the conformal group of $N_{0}$,
$P_{0}=L_{0}N_{0}$ for the Siegel parabolic subgroup of $G_{0}$, etc.

Below are the examples of $N_{0}$ and $G_{0}$ for several different groups $G:$

\begin{itemize}
\item  For $G=O_{p+2,p+2}$, we have $N_{0}=\mathbb{R}$ (rank 1 Jordan
algebra), and $G_{0}=GL_{2}(\mathbb{R)}$.

\item  If $G=Sp_{n,n}$ , then $G_{0}=Sp_{n-1,n-1}$.

\item  If $G=E_{7(7)}$, then $N_{0}=\mathbb{R}^{6,6}$ (rank 2 Jordan algebra),
and $G_{0}=O_{6,6}.$

\item  If $G=E_{7}\frak{(\mathbb{C})}$, then $G_{0}=O_{12}(\mathbb{C)}$.
\end{itemize}

Set $\frak{f}=\bigoplus_{i=2}^{n}\frak{g}^{\varepsilon_{1}-\varepsilon_{i}%
}\oplus\bigoplus_{i=2}^{n}\frak{g}^{\varepsilon_{1}+\varepsilon_{i}}$ and
$\frak{n}^{\prime}=\frak{f}+\frak{n}_{1}$. Then $\frak{n}^{\prime}$ is a
two-step nilpotent subalgebra of $\frak{g}$ with the center $\frak{n}_{1}$.

Any generic unitary irreducible representation of the group $N^{\prime}$ is
determined by the unitary character of its center $N_{1}.$ We denote by
$\rho_{t}$ the unitary irreducible representation of $N^{\prime}$ which
restricts to the multiple of the character $\chi_{t}\,$on $N_{1}$, $t\in
N_{1}^{\vee}=N_{1}^{\ast}\backslash\{0\}$.

Consider now a subgroup $G_{0}N^{\prime}$ of $G$. We can view $G_{0}$ as a
subgroup of a symplectic group $Sp(\frak{f})$ associated with the standard
skew-symmetric bilinear form on $\frak{f}$. Hence we can use the oscillator
representation of $Sp(\frak{f})$ to extend the representation $\rho_{t}$ of
$N^{\prime}$ to a representation of $G_{0}N^{\prime}$ which we denote by
$\widetilde{\rho}_{t}$.

Let $\sigma$ be a unitary representation of $G$, $\operatorname*{rank}%
\nolimits_{N}\sigma=r$, $0<r<n$. Without loss of generality we may assume that
$\sigma$ has no $N_{1}$-fixed vectors. Then by Mackey theory, we can write
down the decomposition%
\[
\sigma|_{G_{0}N^{\prime}}=\int_{N_{1}^{\vee}}^{\oplus}\kappa_{t}%
\otimes\widetilde{\rho}_{t}dt\text{, }%
\]
where all $\kappa_{t}$ are unitary representations of $G_{0}$.

Proceeding as in \cite[3.1]{tens}, we verify that all of the representations
$\kappa_{t}$ are in turn the \emph{low-rank representations of }$G_{0}$. More
precisely, we have the following

\begin{lemma}
\label{=lemma-onedown}Let $\sigma$ be a low-rank representation of $G$,
$\operatorname*{rank}\nolimits_{N}\sigma=r$, $0<r<n$. Then for any $t\in
N_{1}^{\vee}\,$the $N_{0}$-spectrum of the representation $\kappa_{t}$ is
supported on a single $L_{0}$-orbit, and $\operatorname*{rank}\nolimits_{N_{0}%
}\kappa_{t}=r-1$.
\end{lemma}

The next technical lemma is proved exactly as in \cite[3.2]{tens}:

\begin{lemma}
If for all $t\in N_{1}^{\vee}\,$one has $\mathcal{A}(\kappa_{t},G_{0}%
)=\mathcal{A}(\kappa_{t},P_{0})$, then%
\[
\mathcal{A}\left(  \int_{N_{1}^{\vee}}^{\oplus}\kappa_{t}\otimes
\widetilde{\rho}_{t}dt,G_{0}\right)  \subseteq\mathcal{A}\left(  \int
_{N_{1}^{\vee}}^{\oplus}\kappa_{t}\otimes\widetilde{\rho}_{t}dt,P_{0}%
N^{\prime}\right)  \text{.}%
\]

\begin{theorem}
Let $\sigma$ be a representation of $G$, $\operatorname*{rank}\nolimits_{N}%
\sigma=r$, $0<r<n$. Then%
\[
\mathcal{A}(\sigma,G)=\mathcal{A}(\sigma,P)\text{.}%
\]
\end{theorem}
\end{lemma}

\begin{proof}
The proof of the theorem is based on the fact that $G_{0}$ and $P$ generate
$G$, so it is enough to verify that $\mathcal{A}(\sigma,G_{0})\subseteq
\mathcal{A}(\sigma,P)$. Since $P_{0}N^{\prime}$ is a subgroup of $P,$ by the
Lemma above the assertion of the theorem is equivalent to the claim that for
any $t\in N_{1}^{\vee}$
\[
\mathcal{A}(\kappa_{t},G_{0})=\mathcal{A}(\kappa_{t},P_{0})\text{.}%
\]
By the Lemma \ref{=lemma-onedown}, all $\kappa_{t}$ have rank $r-1$.
Proceeding in the same manner, we reduce the statement of the theorem to that
about rank $0$ representations of the certain (classical) group $G_{00}$.
Since all rank 0 representations are the direct integrals of characters, and
any character of $G_{00}$ is determined by its restriction to the Siegel
parabolic $P_{00}\subset G_{00}$, the theorem follows.
\end{proof}

We now consider the tensor product%
\[
\Pi=\pi_{\mathcal{O}^{1}}\otimes\ldots\otimes\pi_{\mathcal{O}^{s}}\text{ ,}%
\]
for $k=k_{1}+k_{2}+\cdots+k_{s}<n$. Then $\Pi$ is low-rank representation of
$G$ and $\operatorname*{rank}\nolimits_{N}\Pi=k$. Applying the theorem above
to $\Pi$, we obtain the statement of Proposition \ref{=prop-vN} for $k<n$.

It remains to check Proposition \ref{=prop-vN} for $k=n$. For all groups $G$
except $O_{p+2,p+2}$, $O_{p+4}(\mathbb{C})$ and $E_{7(7)},$ $E_{7}%
(\mathbb{C})$ the statement of the proposition follows from the results of
\cite{li-lowrank}, since all the representations form the spectrum of $\Pi$
appear in the Howe duality correspondence for appropriate stable range dual
pairs $(G^{1},G)$. For the exceptional cases listed above, the argument can be
constructed along the lines of Section 4 of \cite{tens}.\label{app-lowrank}

\subsection{Tables of groups and symmetric spaces}

In the table below we list the various groups $G$ arising from the
Tits-Kantor-Koecher construction for non-Euclidean Jordan algebras; together
with the associated symmetric spaces $K/M$ and $L/H$; and the crucial root
multiplicities $d$ and $e$ in $\Sigma(\frak{t}_{\mathbb{C}},\frak{k}%
_{\mathbb{C}})$. The rank of the Jordan algebra is $n$ in each case, except
for $O_{p+4}(\mathbb{C})$, $O_{p+2,p+2}$ where the rank is 2, and $E_{7(7)}$,
$E_{7}(\mathbb{C})$ where the rank is 3.%

\[%
\begin{tabular}
[c]{|c|c|c|c|c|}\hline
$G$ & $K/M$ & $L/H$ & $d$ & $e$\\\hline\hline
\multicolumn{1}{|l|}{$GL_{2n}(\mathbb{R})$} & \multicolumn{1}{|l|}{$O_{2n}%
/(O_{n}\times O_{n})$} & \multicolumn{1}{|l|}{$GL_{n}(\mathbb{R)\times}%
GL_{n}(\mathbb{R)}/GL_{n}(\mathbb{R)}$} & \multicolumn{1}{|l|}{$1$} &
\multicolumn{1}{|l|}{$0$}\\\hline
\multicolumn{1}{|l|}{$O_{2n,2n}$} & \multicolumn{1}{|l|}{$(O_{2n}\times
O_{2n})/O_{2n}$} & \multicolumn{1}{|l|}{$GL_{2n}(\mathbb{R)}/Sp_{n}%
(\mathbb{R})$} & \multicolumn{1}{|l|}{$2$} & \multicolumn{1}{|l|}{$0$}\\\hline
\multicolumn{1}{|l|}{$E_{7(7)}$} & \multicolumn{1}{|l|}{$SU_{8}/Sp_{4}$} &
\multicolumn{1}{|l|}{$\mathbb{R}^{\ast}\mathbb{\times}E_{6(6)}/F_{4(4)}$} &
\multicolumn{1}{|l|}{$4$} & \multicolumn{1}{|l|}{$0$}\\\hline
\multicolumn{1}{|l|}{$O_{p+2,p+2}$} & \multicolumn{1}{|l|}{$[O_{p+2}%
]^{2}/\mathbb{[}O_{1}\times O_{p+1}^{2}]$} & \multicolumn{1}{|l|}{$\mathbb{R}%
^{\ast}\times O_{p+1,p+1}/O_{p,p+1}$} & \multicolumn{1}{|l|}{$p$} &
\multicolumn{1}{|l|}{$0$}\\\hline\hline
\multicolumn{1}{|l|}{$Sp_{n}(\mathbb{C})$} & \multicolumn{1}{|l|}{$Sp_{n}%
/U_{n}$} & \multicolumn{1}{|l|}{$GL_{n}(\mathbb{C)}/O_{n}(\mathbb{C})$} &
\multicolumn{1}{|l|}{$1$} & \multicolumn{1}{|l|}{$1$}\\\hline
\multicolumn{1}{|l|}{$GL_{2n}(\mathbb{C})$} & \multicolumn{1}{|l|}{$U_{2n}%
/(U_{n}\times U_{n})$} & \multicolumn{1}{|l|}{$GL_{n}(\mathbb{C)\times}%
GL_{n}(\mathbb{C)}/GL_{n}(\mathbb{C)}$} & \multicolumn{1}{|l|}{$2$} &
\multicolumn{1}{|l|}{$1$}\\\hline
\multicolumn{1}{|l|}{$O_{4n}(\mathbb{C})$} & \multicolumn{1}{|l|}{$O_{4n}%
/U_{2n}$} & \multicolumn{1}{|l|}{$GL_{2n}(\mathbb{C)}/Sp_{n}(\mathbb{C})$} &
\multicolumn{1}{|l|}{$4$} & \multicolumn{1}{|l|}{$1$}\\\hline
\multicolumn{1}{|l|}{$E_{7}(\mathbb{C})$} & \multicolumn{1}{|l|}{$E_{7}%
/(E_{6}\times U_{1})$} & \multicolumn{1}{|l|}{$\mathbb{C}^{\ast}%
\mathbb{\times}E_{6}(\mathbb{C})/F_{4}(\mathbb{C})$} &
\multicolumn{1}{|l|}{$8$} & \multicolumn{1}{|l|}{$1$}\\\hline
\multicolumn{1}{|l|}{$O_{p+4}(\mathbb{C})$} & \multicolumn{1}{|l|}{$O_{p+4}%
/(O_{p+2}\times U_{1})$} & \multicolumn{1}{|l|}{$\mathbb{C}^{\ast}\times
O_{p+2}(\mathbb{C})/O_{p+1}(\mathbb{C})$} & \multicolumn{1}{|l|}{$p$} &
\multicolumn{1}{|l|}{$1$}\\\hline\hline
\multicolumn{1}{|l|}{$Sp_{n,n}$} & \multicolumn{1}{|l|}{$(Sp_{n}\times
Sp_{n})/Sp_{n}$} & \multicolumn{1}{|l|}{$GL_{n}(\mathbb{H)}/O_{n}^{\ast}$} &
\multicolumn{1}{|l|}{$2$} & \multicolumn{1}{|l|}{$2$}\\\hline
\multicolumn{1}{|l|}{$GL_{2n}(\mathbb{H)}$} & \multicolumn{1}{|l|}{$Sp_{2n}%
/(Sp_{n}\times Sp_{n})$} & \multicolumn{1}{|l|}{$GL_{n}(\mathbb{H)\times
}GL_{n}(\mathbb{H)}/GL_{n}(\mathbb{H)}$} & \multicolumn{1}{|l|}{$4$} &
\multicolumn{1}{|l|}{$3$}\\\hline
\end{tabular}
\]
In this table we list the various homogeneous spaces $X=G^{\prime}/H^{\prime}$
that arise in the $\theta$-correspondence of Theorem \ref{theoremB}.
\[%
\begin{tabular}
[c]{|c|c|}\hline
$G$ & $X$\\\hline\hline
\multicolumn{1}{|l|}{$GL_{2n}(\mathbb{R})$} & \multicolumn{1}{|l|}{$GL_{k}%
(\mathbb{R)}/\left[  GL_{k_{1}}(\mathbb{R)}\times\cdots\times GL_{k_{s}%
}(\mathbb{R)}\right]  $}\\\hline
\multicolumn{1}{|l|}{$O_{2n,2n}$} & \multicolumn{1}{|l|}{$Sp_{k}%
(\mathbb{R})/\left[  Sp_{k_{1}}(\mathbb{R)}\times\cdots\times Sp_{k_{s}%
}(\mathbb{R)}\right]  $}\\\hline
\multicolumn{1}{|l|}{$E_{7(7)}$} & \multicolumn{1}{|l|}{$%
\begin{array}
[c]{lc}%
Spin_{4,5}/Spin_{4,4} & (k_{1}=1,k_{2}=1)\\
F_{4(4)}/Spin_{4,5} & (k_{1}=2,k_{2}=1)
\end{array}
$}\\\hline
\multicolumn{1}{|l|}{$O_{p+2,p+2}$} & \multicolumn{1}{|l|}{$SO_{p,p+1}%
/SO_{p,p}$ $(k_{1}=1,k_{2}=1)$}\\\hline\hline
\multicolumn{1}{|l|}{$Sp_{n}(\mathbb{C})$} & \multicolumn{1}{|l|}{$O_{k}%
(\mathbb{C})/\left[  O_{k_{1}}(\mathbb{C})\times\cdots\times O_{k_{s}%
}(\mathbb{C})\right]  $}\\\hline
\multicolumn{1}{|l|}{$GL_{2n}(\mathbb{C})$} & \multicolumn{1}{|l|}{$GL_{k}%
(\mathbb{C)}/\left[  GL_{k_{1}}(\mathbb{C)}\times\cdots\times GL_{k_{s}%
}(\mathbb{C)}\right]  $}\\\hline
\multicolumn{1}{|l|}{$O_{4n}(\mathbb{C})$} & \multicolumn{1}{|l|}{$Sp_{k}%
(\mathbb{C})/\left[  Sp_{k_{1}}(\mathbb{C)}\times\cdots\times Sp_{k_{s}%
}(\mathbb{C)}\right]  $}\\\hline
\multicolumn{1}{|l|}{$E_{7}(\mathbb{C})$} & \multicolumn{1}{|l|}{$\mathbb{\ }%
\begin{array}
[c]{lc}%
Spin_{9}(\mathbb{C})/Spin_{8}(\mathbb{C}) & (k_{1}=1,k_{2}=1)\\
F_{4}(\mathbb{C})/Spin_{9}(\mathbb{C}) & (k_{1}=2,k_{2}=1)
\end{array}
$}\\\hline
\multicolumn{1}{|l|}{$O_{p+4}(\mathbb{C})$} & \multicolumn{1}{|l|}{$SO_{p+1}%
(\mathbb{C})/SO_{p}(\mathbb{C})$ $(k_{1}=1,k_{2}=1)$}\\\hline\hline
\multicolumn{1}{|l|}{$Sp_{n,n}$} & \multicolumn{1}{|l|}{$O_{k}^{\ast}/\left[
O_{k_{1}}^{\ast}\times\cdots\times O_{k_{s}}^{\ast}\right]  $}\\\hline
\multicolumn{1}{|l|}{$GL_{2n}(\mathbb{H)}$} & \multicolumn{1}{|l|}{$GL_{k}%
(\mathbb{H)}/\left[  GL_{k_{1}}(\mathbb{H)}\times\cdots\times GL_{k_{s}%
}(\mathbb{H)}\right]  $}\\\hline
\end{tabular}
\]

\label{sec-table}


\begin{thebibliography}{DS1}
\bibitem[A]{adams}Adams, J.F., \emph{Lectures on exceptional Lie groups,
}University of Chicago Press, Chicago 1996

\bibitem[BK]{braun}Braun, H. and Koecher, M., \emph{Jordan-Algebren},
Springer, Berlin -- New York 1966

\bibitem[DS1]{tens}Dvorsky, A. and Sahi, S.,\emph{\ Tensor products of
singular representations and an extension of the }$\theta$%
\emph{-correspondence, }Selecta Math. \textbf{4 }(1998), 11-29

\bibitem[DS2]{hilbert-two}Dvorsky, A. and Sahi, S., \emph{Explicit Hilbert
spaces for certain unipotent representations II}, Invent. Math. \textbf{138
}(1999), 203--224

\bibitem[H]{helgason}Helgason, S., \emph{Groups and Geometric Analysis},
Mathematical surveys and monographs v. 83,\textbf{\ }Amer. Math. Soc.,
Providence 2000

\bibitem[KS]{kostant-sahi}Kostant, B. and Sahi, S., \emph{Jordan algebras and
Capelli identities,} Invent. Math. \textbf{112} (1993), 657--664

\bibitem[Li1]{li-singular}Li, J.-S., \emph{Singular unitary representations of
classical groups, }Invent. Math. \textbf{97 }(1989), 237--255

\bibitem[Li2]{li-lowrank}Li, J.-S., \emph{On the classification of irreducible
low rank unitary representations of classical groups, }Comp. Math \textbf{71
}(1989), 29-48

\bibitem[Lo]{loos}Loos, O., \emph{Bounded symmetric domains and Jordan pairs},
Mathematical Lectures, University of California, Irvine 1977

\bibitem[M]{mackey}Mackey, G., \emph{The theory of unitary group
representations}, University of Chicago Press, Chicago 1976

\bibitem[Op]{opdam}Opdam, E., \emph{Dunkl operators, Bessel functions and the
discriminant of a finite Coxeter group}, Compositio Math. \textbf{85} (1993), 333-373

\bibitem[OS]{oshima}Oshima, T. and Sekiguchi, J., \emph{The restricted root
system of a semisimple symmetric pair}, In: Adv. Stud. Pure Math. \textbf{4},
433-497, North-Holland, Amsterdam 1984

\bibitem[P]{poguntke}Poguntke, D., \emph{Unitary representations of Lie groups
and operators of finite rank}, Ann. of Math. (2) \textbf{140} (1994), 503-556

\bibitem[S1]{sahi-expl}Sahi, S., \emph{Explicit Hilbert spaces for certain
unipotent representations}, Invent. Math. \textbf{110 }(1992), 409--418

\bibitem[S2]{shilov}Sahi, S., \emph{Unitary representations on the Shilov
boundary of a symmetric tube domain}, In: Contemp. Math. \textbf{v. 145}
\textbf{\ }(1993), 275-286

\bibitem[S3]{sahi-dp}Sahi, S., \emph{Jordan algebras and degenerate principal
series}, J. reine angew. Math. \textbf{462} (1995), 1--18

\bibitem[SS]{sahi-stein}Sahi, S. and Stein, E., \emph{Analysis in matrix space
and Speh's representation}, Invent. Math. \textbf{101 }(1990), 379-393

\bibitem[Sc]{sch}Schlichtkrull, H., \emph{Hyperfunctions and harmonic analysis
on symmetric spaces}, Birkh\"{a}user, Boston 1984

\bibitem[Sh]{shimura}Shimura, G., \emph{Generalized Bessel functions on
symmetric spaces}, J. reine angew. Math. \textbf{509} (1999), 35-66

\bibitem[V]{vogan-nogo}Vogan, D., \emph{Singular unitary representations,} In:
Lect. Notes Math. \textbf{880}, 506-535, Springer, Berlin-New York 1980
\end{thebibliography}
\end{document}